\scrollmode
\documentclass[12pt]{amsart}
\usepackage{verbatim}
\usepackage{eucal}
\usepackage{amssymb}
\usepackage{graphicx}
\usepackage{psfrag}
\usepackage{mathrsfs}
\usepackage[all]{xy}
\newif \ifwide
\newif \ifavnermargin
\def \makemargins{
\ifwide
	\oddsidemargin .25in
	\evensidemargin .25in
	\textwidth 6.00in
\else
\fi
\ifavnermargin
	\headheight=7pt
	\textheight=574pt
	\textwidth=432pt
	\topmargin=14pt
	\oddsidemargin=18pt
	\evensidemargin=18pt
\else	
\fi
}

\avnermargintrue
\makemargins
\swapnumbers
\theoremstyle{plain}
\newtheorem{theorem}[subsection]{Theorem}
\newtheorem{conjecture}[subsection]{Conjecture}
\newtheorem{proposition}[subsection]{Proposition}
\newtheorem{lemma}[subsection]{Lemma}

\theoremstyle{definition}
\newtheorem{definition}[subsection]{Definition}
\newtheorem{algorithm}[subsection]{Algorithm}

\theoremstyle{remark}
\newtheorem{remark}[subsection]{Remark}

\newcount\timehh\newcount\timemm
\timehh=\time 
\divide\timehh by 60 \timemm=\time
\newcommand{\draftauthor}[1]{\author{#1
    {
      --- \protect \protect\sc\today\ ---
      \ifnum\timehh<10 0\fi\number\timehh\,:\,\ifnum\timemm<10 0\fi\number\timemm
      \protect \, \, \protect \bf DRAFT
    }
  }
}
\newcommand{\R}{{\mathbb R}}

\newcommand{\Z}{{\mathbb Z}}
\newcommand{\Q}{{\mathbb Q}}

\newcommand{\V}{{\mathscr V}}

\newcommand{\MMM}{{\mathscr M}}
\newcommand{\HHH}{{\mathscr H}}

\newcommand{\PPP}{{\mathscr P}}

\newcommand{\OK}{{\mathscr O_{K}}}

\newcommand{\Vor}{Vorono\v{\i}}

\newcommand{\vv}{{\mathbf {v}}}

\newcommand{\uu}{{\mathbf {u}}}

\newcommand{\coinv}[1]{(S_{#1})_{\Gamma }}

\DeclareMathOperator*{\Max}{Max}
\DeclareMathOperator*{\Min}{Min}
\DeclareMathOperator{\support}{supp}

\DeclareMathOperator{\Cand}{cand}
\DeclareMathOperator{\Vertex}{vert}

\DeclareMathOperator{\diag}{diag}
\DeclareMathOperator{\sign}{sgn}
\DeclareMathOperator{\vol}{vol}
\CompileMatrices

\begin{document}

\title[Computing Hecke eigenvalues]{Computing Hecke eigenvalues below the cohomological dimension}

\newif \ifdraft
\def \makeauthor{
\ifdraft
	\draftauthor{Paul E. Gunnells}
\else
\author{Paul E. Gunnells}
\address{Department of Mathematics\\
Columbia University\\
New York, NY  10027}
\email{gunnells@math.columbia.edu}
\fi
}

\thanks{The author was partially supported by a Columbia University Faculty
Research grant and NSF Grant DMS 96--27870.}
\draftfalse
\makeauthor

\ifdraft
	\date{\today}
\else
	\date{November, 1998.  Revised September, 1999}
\fi

%
%
\begin{abstract}
Let $\Gamma$ be a torsion-free finite-index subgroup of $SL_{n} (\Z )$
or $GL_{n} (\Z )$, and let $\nu $ be the cohomological dimension of
$\Gamma $.  We present an algorithm to compute the eigenvalues of the
Hecke operators on $H^{\nu -1} (\Gamma ;\Z )$, for
$n=2$,~$3$,~and~$4$.  In addition, we describe a modification of the
modular symbol algorithm of Ash-Rudolph \cite{ash.rudolph} for
computing Hecke eigenvalues on $H^{\nu } (\Gamma ;\Z )$.
\end{abstract}
\maketitle

%
%
\section{Introduction}\label{introduction}
\subsection{}
Let $\Gamma$ be a finite-index subgroup of $SL_{n} (\Z )$ or $GL_{n}
(\Z )$, and let $\MMM$ be a $\Z \Gamma $-module.  The group cohomology
$H^{*} (\Gamma ;\MMM )$ plays an important role in number theory,
through its connection with automorphic forms and representations of
the absolute Galois group.  For an introduction to this conjectural
framework, see \cite{ash:galois}.

For $n=2$ and $\Gamma $ a congruence subgroup, the arithmetic nature
of $H^{*} (\Gamma; \MMM )$ has been decisively confirmed
(cf. \cite{shimura}).  For higher dimensions the picture is
mysterious, although several compelling examples for $n=3$ have
appeared recently in the literature.  In \cite{apt}, rational
cohomology classes of certain $\Gamma \subset GL_{3}
(\Z )$ are related to modular Galois representations.
Many more examples of this phenomenon appear in \cite{exp.ind}.  In
\cite{geemen.top, vgt2}, rational cohomology classes of certain congruence
groups are related to the Hasse-Weil zeta functions of certain
surfaces.  Finally, in \cite{aac} torsion classes in the cohomology of
$\Gamma = GL_{3} (\Z )$ with twisted coefficients are linked to
modular Galois representations, and in \cite{ash.tiep} the arithmetic
nature of many of these classes is proven.

\subsection{}
In all cases, the arithmetic significance of $H^{*} (\Gamma; \MMM )$
is revealed through the \emph{Hecke operators}.  These are
endomorphisms of the cohomology associated to certain finite-index
subgroups of $\Gamma $.  The eigenvalues of these linear maps provide
a ``signature'' for the cohomology, which one hopes can be matched to
number-theoretic data.  Thus to test these conjectures, or to search
for counterexamples, it is crucially important to be able to compute
Hecke eigenvalues.

\subsection{}\label{intro2}
In general, computing these eigenvalues is a difficult problem.
Essentially the only technique available in the literature is
the \emph{modular symbol algorithm} (\S\ref{mod.symb.section}), due to
Manin \cite{manin} ($n=2$) and Ash and Rudolph \cite{ash.rudolph}
($n\geq 3$).  Using this algorithm one can compute the Hecke action on
$H^{\nu } (\Gamma; \MMM )$, where $\nu $ is the \emph{cohomological
dimension} of $\Gamma $.  That is, $\nu $ is the smallest number such
that $H^{i} (\Gamma ; \MMM ) =0$ for $i>\nu $ and any $\MMM $.

In particular, if $\Gamma \subset SL_{3} (\Z )$ or $GL_{3} (\Z )$,
then $\nu =3$.  This is the focus of \cite{apt, exp.ind, geemen.top,
vgt2, aac, ash.tiep}.  For certain congruence groups $\Gamma$, the
groups $H^{3} (\Gamma; \Q )$ and $H^{2}(\Gamma ;\Q )$ contain
\emph{cuspidal} cohomology classes \cite{agg}.  In a certain sense
these classes are the most interesting constituents of the cohomology.
A Lefschetz duality argument \cite[Theorem 3.1]{ash.tiep} shows that
the cuspidal eigenclasses in $H^{2} (\Gamma; \Q )$ have the same
eigenvalues as those in $H^{3} (\Gamma;\Q )$, and therefore the
modular symbol algorithm suffices to compute Hecke eigenvalues in this
dimension.

Now suppose $\Gamma \subset SL_{4} (\Z )$, so that $\nu =6$.  In this
case $H^{6} (\Gamma; \Q )$ does not contain cuspidal classes, and one
is interested in $H^{5} (\Gamma;\Q )$.  Again one wants to compute
Hecke eigenvalues, and since Lefschetz duality relates $H^{6}$ to
$H^{3}$ and not $H^{5}$, the modular symbol algorithm doesn't apply.
Thus one has the natural problem of devising an algorithm for this
context.

\subsection{}
The purpose of this article is to describe an algorithm that---in
practice---allows computation of the Hecke action on $H^{\nu -1}
(\Gamma ;\Z )$, where $\Gamma $ is a torsion-free subgroup of $SL_{n}
(\Z )$ or $GL_{n} (\Z )$, and $n\leq 4$.  We emphasize that the phrase
``in practice'' is to be taken literally.

Let us be more precise.  To represent elements of $H^{\nu -1} (\Gamma
;\Z )$, we use chains in the \emph{sharbly complex}\footnote{The name
is due to Lee Rudolph, in honor of the authors of \cite{sharblee1}.}
$S_{*}$ (\S\ref{sharb.stuff}).  This is a complex of infinite $\Gamma
$-modules such that the homology of the complex of coinvariants
$\coinv{*}$ is naturally isomorphic to $H^{*} (\Gamma ;\Z )$.
Furthermore, $S_{*}$ has a natural Hecke action (\S\ref{relationship})
that passes to $\coinv{*}$.  Hence the sharbly complex provides a
convenient setting to study the cohomology as a Hecke module.

Both $S_{*}$ and $\coinv{*}$ are infinitely generated, whereas
$H^{*} (\Gamma ; \Z )$ is finitely generated.  Hence for practical
computations we must identify a finite subset of $\coinv{*}$ that
spans the cohomology.  For $n\leq 4$ and $H^{\nu -1}$, a spanning set
is provided by the \emph{reduced sharblies} (\S\ref{one.sharb}), which
form a subcomplex of $S_{*}$.  Unfortunately, the Hecke operators do
not preserve this subcomplex.  Thus to compute eigenvalues, we must
describe an algorithm that writes a general sharbly cycle as a sum of
reduced sharbly cycles.

\subsection{}
So suppose $\xi $ is a sharbly cycle mod $\Gamma $ representing a
class in $H^{\nu -1} (\Gamma ;\Z )$.  There is a function
$\|\phantom{\uu }\|\colon S_{*}\rightarrow \Z $ such that $\xi $ is
reduced if and only if $\|\xi\|=1$ (Definition \ref{reduced.sharbly}).
Algorithm~\ref{hecke.algorithm} describes a process that takes $\xi$
as input and produces a cycle $\xi '$ homologous to $\xi $ in
$\coinv{*}$.  Geometrically, the algorithm acts by applying the
modular symbol algorithm simultaneously over all of $\xi $.  Of
course, to be useful for eigenvalue computations, we want that if
$\|\xi\|>1$, then $\|\xi '\|<\|\xi \|$.

We cannot prove that the output $\xi '$ will satisfy this inequality.
However, for $n\leq 4$---the cases of practical interest---this
inequality has always held.  More precisely, in computer experiments
(\S\ref{experiments.section}) with both random data and $1$-sharbly
cycles for $n\leq 4$, Algorithm \ref{hecke.algorithm} has always
successfully written a general $1$-sharbly cycle as a sum of reduced
$1$-sharbly cycles.  Currently we are applying Algorithm
\ref{hecke.algorithm} in joint work with Avner Ash and Mark McConnell
to decompose $H^{5} (\Gamma ;\Q )$ as a Hecke module for certain
congruence groups $\Gamma \subset SL_{4} (\Z)$  \cite{agm}.  Details of these
computations will appear in a later publication.

\subsection{} 
Here is a guide to this paper.  In \S\ref{background} we recall the
topological and combinatorial background necessary for computing
$H^{*} (\Gamma ;\Z )$.  We discuss the reduction theory due to \Vor \
\cite{voronoi1} and the sharbly complex, as well as the Hecke
operators and how they interact with the sharbly complex.  In
\S\ref{modular.symbols} we recall the modular symbol algorithm, and
describe two new conjectural techniques to implement it (Conjectures
\ref{vorconj} and~\ref{lllconj}).  These techniques link the modular
symbol algorithm to \Vor \ reduction and $LLL$-reduction, and are
conjectured to be true in all dimensions.  We also include proofs of
the conjectures in special cases.  Then in \S\ref{one.sharb} we
present Algorithm~\ref{hecke.algorithm} and prove that, given a
sharbly cycle $\xi $ mod $\Gamma $ as input, the output $\xi '$ is a
homologous cycle mod $\Gamma $ (Theorem~\ref{main.thm}).  We also
discuss conditions under which we expect $\|\xi '\|<\|\xi \|$
(Conjecture~\ref{os.conj}).  Finally, in \S\ref{experiments.section},
we describe experiments we performed to generate evidence for
Conjectures~\ref{vorconj},~\ref{lllconj}, and~\ref{os.conj}.

\subsection{Acknowledgements}\label{thanx} 
We thank Romuald Dabrowski and Lee Rudolph for helpful discussions.
We are especially grateful to Avner Ash and Mark McConnell for much
enthusiastic help and advice during this project.  We thank the
referee, whose many comments vastly improved the exposition.  Finally,
we thank Bob MacPherson for originally suggesting the problem and for
much encouragement and support.

The experiments performed in this paper to develop Algorithm
\ref{hecke.algorithm} and to test Conjectures~\ref{vorconj},
\ref{lllconj}, and~\ref{os.conj} were implemented using several
software packages: GP-Pari~\cite{pari}, LiDIA~\cite{LiDIA},
Mathematica~\cite{mathematica}, and SHEAFHOM~\cite{sheafhom}.  

%
%
%
\section{Background}\label{background}
In this section we describe the topological tools we use to study
$H^{*} (\Gamma ; \Z )$: the \Vor \ polyhedron and the sharbly complex.
We present these objects in the context of $\Gamma \subset SL_{n} (\Z
)$.  However, all of what we say applies with minor modification to
$\Gamma \subset GL_{n} (\Z )$.

\subsection{}\label{vec.space}
Let $V$ be the $\R $-vector space of symmetric $n\times n$
matrices, and let $C\subset V$ be the cone of positive-definite
matrices.  The linear group $G=SL_{n} (\R )$ acts on $C$ by
$(g,c)\mapsto g\cdot c\cdot g^{t}$, and the stabilizer of any given
point is isomorphic to $SO_{n}$.

Let $X$ be $C$ mod homotheties.  The $G$-action on $C$ commutes
with the homotheties and induces a transitive $G$-action on $X$.  The
stabilizer of any fixed point of $X$ is again $SO_{n}$.  After
choosing a basepoint, we may identify $X$ with the global Riemannian
symmetric space $SL_{n} (\R )/SO_{n}$, a contractible, noncompact,
smooth manifold of real dimension $N = n (n+1)/2 - 1$.

The group $SL_{n} (\Z )$ acts on $X$ via the $G$-action, and does
so properly discontinuously.  Hence if $\Gamma \subset SL_{n} (\Z )$
is any torsion-free subgroup, the quotient $\Gamma \backslash X$ is
a real noncompact manifold, and is an Eilenberg-Mac~Lane space for
$\Gamma $.  We may then identify the group cohomology $H^{*} (\Gamma;
\Z )$ with $H^{*} (\Gamma \backslash X;\Z )$.  Although the dimension
of $\Gamma \backslash X$ is $N$, it can be shown that $H^{i} (\Gamma
\backslash X;\Z )=0$ if $i\geq N-n+1$ \cite[Theorem
11.4.4]{borel.serre}.  The number $\nu = N-n+1$ is called the
\emph{cohomological dimension} of $\Gamma $.

\subsection{}\label{voronoi}
Recall that a point in $\Z ^{n}$ is said to be \emph{primitive} if the
greatest common divisor of its coordinates is $1$.  In particular, a
primitive point is nonzero.  Let $\PPP\subset \Z ^{n}$ be the
set of primitive points.  Any $v\in \PPP$, written as a column vector,
determines a rank-one quadratic form $q (v)\in \bar C$ by $q ( v) =
v\cdot v^{t}$.
\begin{definition}\label{vorpoly}
The \emph{\Vor \ polyhedron} $\Pi $ is the closed convex hull of the
points $q (v)$, as $v$ ranges over $\PPP$.
\end{definition}
Note that, by construction, $SL_{n} (\Z )$ acts on $\Pi $.  The cones
over the faces of $\Pi $ form a fan $\V$ that induces a $\Gamma
$-admissible decomposition of $C$~\cite[p. 117]{ash}.  Essentially,
this means that $\Gamma $ acts on $\V$; that each cone is spanned by a
\emph{finite} collection of points $q (v)$ where $v\in \PPP$; and that
mod $\Gamma $ there are only finitely many orbits in $\V$.  The fan
$\V$ provides a reduction theory for $C$ in the following sense: any
point $x\in C$ is contained in a unique $\sigma \in \V$.

Given $\sigma \in \V$, let $\Vertex \sigma $ be the set of all $v\in
\PPP$ such that $q (v)$ is a vertex of the face of $\Pi $ generating
$\sigma $.  For later use, we record the following theorem of \Vor :

\begin{theorem}\label{vor.thm}
\cite{voronoi1} Let $E$ be the standard basis of $\Z ^{n}$, and let
$\Sigma$ be the cone spanned by the $n (n+1)/2$ points $q (e_{i})$ and
$q (e_{i}-e_{j})$, where $e_{i}\in E$  and $1\leq i<j\leq n$.  Then
$\Sigma$ occurs as a top-dimensional cone in $\V $ for all $n$.
\end{theorem}

\subsection{}\label{sharb.stuff}
We now discuss an algebraic tool to compute $H^{*} (\Gamma ;\Z )$.  The
material in this section closely follows \cite{ash.sharb}.

Recall that the \emph{Steinberg module} $St(n)$ is the $\Z \Gamma
$-module $H^{\nu }(\Gamma ;\Z \Gamma )$.  
\begin{theorem}\label{msymb}
\cite{ash.sharb} The Steinberg module is isomorphic to the module of
formal $\Z $-linear combinations of the elements $[v_{1},
\ldots, v_{n}]^{*}$, where each $v_{i}\in \Q ^{n}$ is nonzero, mod
the relations:
\begin{enumerate}\label{mod.sym.def}
\item If $\tau$ is a permutation on $n$ letters, then
$[v_{1},\ldots,v_{n}]^{*} = \sign (\tau ) [\tau (v_{1}),\ldots,\tau
(v_{n})]^{*}$, where $\sign (\tau)$ is the sign of $\tau $.
\item If $q\in \Q ^{\times }$, then $[q v_{1},v_{2},\ldots,v_{n}]^{*}
= [v_{1},\ldots,v_{n}]^{*}$.
\item If the $v_{i}$ are linearly dependent, then
$[v_{1},\ldots,v_{n}]^{*} = 0$.
\item If $v_{0},\ldots,v_{n}$ are nonzero points in $\Q^{n} $, then
$\sum _{i}(-1)^{i}[v_{0},\ldots,\hat{v_{i}},\ldots,v_{n}]^{*}= 0$.
\end{enumerate}
\end{theorem}

By Borel-Serre duality \cite[\S11.4]{borel.serre}, if $\Gamma $ is
torsion-free, then for any $\Z \Gamma $-module $\MMM $ we have a natural
isomorphism
\begin{equation}\label{bs.duality}
\Phi \colon H^{k}(\Gamma ; \MMM ) \longrightarrow H_{\nu - k} (\Gamma
;St(n)\otimes \MMM ).
\end{equation}
Hence one may compute $H^{*} (\Gamma ; \Z )$ by computing the homology
of a $\Z \Gamma$-free resolution of $St (n)\otimes \Z $.  Such a
resolution is provided by the sharbly complex.

\begin{definition}\label{sharbly.complex}
\cite{ash.sharb} The \emph{sharbly complex} is the chain complex
$\left\{S_{*},\partial \right\}$ given by the following data:
\begin{enumerate}
\item For $k\geq 0$, $S_{k}$ is the module of formal $\Z$-linear
combinations of elements $\uu  = [v_{1},\ldots,v_{n+k}]$, where each
$v_{i}\in \PPP$, mod the relations:
\begin{enumerate}
\item If $\tau $ is a permutation on $(n+k)$ letters, then
\[
[v_{1},\ldots,v_{n+k}] = \sign (\tau ) [\tau (v_{1}),\ldots,\tau
(v_{n+k})], 
\]
where $\sign (\tau )$ is the sign of $\tau $.
\item  If $q = \pm 1$, then 
\[
[q v_{1},v_{2}\ldots,v_{n+k}] = [v_{1},\ldots,v_{n+k}].
\]
\item If the rank of the matrix
$(v_{1},\ldots,v_{n+k})$ is less than $n$, then $\uu = 0$.
\end{enumerate}
\item The boundary map $\partial \colon S_{k}\rightarrow S_{k-1}$ is 
$$[v_{1},\ldots,v_{n+k}] \longmapsto \sum _{i=1}^{n+k}  (-1)^{i}
[v_{1},\ldots,\hat{v_{i}},\ldots,v_{n+k}].  $$
\end{enumerate}
\end{definition}

The elements $\uu = [v_{1},\dots ,v_{n+k}]$ are called
\emph{$k$-sharblies}.  A $0$-sharbly is also called a \emph{modular
symbol}.  By abuse of notation, we will often use the same symbol $\uu
$ to denote a $k$-sharbly and the $k$-sharbly chain $1\cdot \uu $.
The obvious left action of $\Gamma $ on $S_{* }$ commutes with
$\partial $.

\begin{proposition}\label{ashsharbprop}
\cite{ash.sharb} The complex $\left\{S_{*},\partial \right\}$ is a a
$\Z \Gamma $-free resolution of $St(n)$, with the map
$S_{0}\rightarrow St(n)$ given by $\uu \mapsto \uu ^{*}$.
\end{proposition}

For any $k\geq 0$, let $\coinv{k}$ be the module of $\Gamma
$-coinvariants.  This is the quotient of $S_{k}$ by the relations of
the form $\gamma \cdot \uu - \uu $, where $\gamma \in \Gamma $, $\uu
\in S_{k}$.  This is also a complex with the induced boundary, which
we denote by $\partial_{\Gamma } $.  Proposition~\ref{ashsharbprop}
and~\eqref{bs.duality} imply that $H^{k}(\Gamma ; \Z)$ is naturally
isomorphic to $H_{\nu -k} (\coinv{*})$.

\subsection{}\label{hecke.operators}
Now we recall the definition of the Hecke operators.  More details can
be found in \cite[Ch. 3]{shimura}.

Fix an arithmetic group $\Gamma \subset SL_{n} (\Z )$.  Given $g\in
GL_{n} (\Q )$, let $\Gamma ^{g}= g^{-1}\Gamma g$ and $\Gamma ' =
\Gamma \cap \Gamma^g$.  Then $[\Gamma : \Gamma ']$ and $[\Gamma^g:
\Gamma']$ are finite.  The inclusions $\Gamma '\rightarrow \Gamma $
and $\Gamma '\rightarrow \Gamma^ g$ determine a diagram
\[
\vbox{\xymatrix{&{\Gamma '\backslash X}\ar[dr]^{t}\ar[dl]_{s}&\\
          {\Gamma \backslash X}&&{\Gamma \backslash X}}}
\]
Here $s(\Gamma 'x) = \Gamma x$ and $t$ is the composition of $\Gamma '
x \mapsto \Gamma ^{g}x$ with left multiplication by $g$.  This diagram
is the {\em Hecke correspondence} associated to $g$.  It can be shown
that, up to isomorphism, the Hecke correspondence depends only on the
double coset $\Gamma g\Gamma $.

Because the maps $s$ and $t$ are proper, they induce a map on
cohomology:
\[
T_{g} := t_{*}s^{*}\colon H^{*}(\Gamma \backslash X;\Z
)\rightarrow H^{*}(\Gamma \backslash X;\Z). 
\]
This is the {\em Hecke operator} associated to $g$.
We let ${\HHH}_{\Gamma }$ be the $\Z $-algebra
generated by the Hecke operators, with product given by composition.

For an example, let $\Gamma = SL_{n} (\Z )$.  Then $\HHH _{\Gamma }$
decomposes as a tensor product
\[
\HHH _{\Gamma } = \bigotimes _{\text{$p$ prime}} \HHH _{p}.
\]
Each $\HHH _{p}$ is a polynomial ring generated by the double cosets
\begin{equation}\label{diag.cosets}
T_{p} (k,n) = \Gamma \diag (1,\dots ,1,\underbrace{p,\dots
,p}_{k})\,\Gamma.
\end{equation}

\subsection{}\label{relationship}
Now let $u\in H^{k} (\Gamma ; \Z )$ be a cohomology class.  Choose $g\in
GL_{n} (\Q )$, and let $T_{g}\in \HHH $ be the Hecke operator
associated to $g$.  We want to explicitly describe the action of
$T_{g}$ on $u$ in terms of the sharbly complex.

Choose $\xi\in S_{k}$ such that $\xi $ is a cycle mod $\Gamma $ and
$\Phi ^{-1} (\xi ) = u$.  Write $\xi = \sum n (\uu ) \uu $, where $n
(\uu )\in \Z $, and almost all $n (\uu )=0$.  The double coset $\Gamma
g\Gamma $ decomposes as
\[
\Gamma  g\Gamma  = \coprod_{h\in I} \Gamma h
\]
for some set $I\subset GL_{n} (\Q )$.  Note that $I$ is finite.  We
have a map $S_{k}\rightarrow S_{k}$ given by
\begin{equation}\label{hecke.map}
T_{g}\colon \xi \longmapsto \sum_{\substack{h\in I}} n (\uu ) h\cdot \uu.
\end{equation}
One can show that the right-hand side of \eqref{hecke.map} is a
well-defined cycle mod $\Gamma $, and that under $\Phi $ this cycle
passes to $T_{g} (u)$.

In general, $I\not \subset SL_{n} (\Z )$.  Thus the Hecke operators do
not preserve the subcomplex of $S_{*}$ generated by $\V $.

%
%
\section{Modular Symbols}\label{modular.symbols}
In this section we recall the Ash-Rudolph modular symbol algorithm and
present our conjectural implementations of it.

\subsection{}\label{mod.symb.section}
Let $\xi $ be a $k$-sharbly chain, and write $\xi =\sum n (\uu )\uu $,
where $n (\uu )\in \Z $ and almost all $n (\uu )=0$.  Let $\support
\xi $ be the set of $k$-sharblies $\{\uu \mid n (\uu )\not = 0 \}$.
Let $Z(\xi )$ be the set of all modular symbols that appear as a
submodular symbol of some $\uu \in \support \xi $.  In other words,
$\vv \in Z (\xi )$ if and only if there is a $\uu = [v_{1},\dots
,v_{n+k}]\in \support \xi $ such that $ \vv =[v_{i_{1}},\dots
,v_{i_{n}}]$ for $\{i_{1},\dots ,i_{n} \}\subset \{1,\dots ,n+k \}$.

\begin{definition}\label{reduced.sharbly}
Given any modular symbol $\vv = [v_{1},\dots ,v_{n}]$, let 
\[
\|\vv\|= |\det (v_{1},\dots ,v_{n})|.
\]
We extend this to
$\|\phantom{\xi } \|\colon S_{k}\rightarrow \Z $ by setting
\[
\|\xi \| = \Max_{\vv \in Z (\xi )}\bigl\{\|\vv \| \bigr\}.
\]
We say $\xi $ is \emph{reduced} if $\|\xi \| = 1$.  In the special
case that $\xi =\vv $ is a modular symbol, we say that $\xi $ is a
\emph{unimodular symbol}.  
\end{definition}

Note that $\|\phantom{\xi } \|$ is well-defined modulo the relations
in Definition~\ref{sharbly.complex}.

The reduced $k$-sharbly chains form a finitely generated subgroup of
$\coinv{k}$.  In general, the image of this subgroup under
the map $S_{k}\rightarrow H^{\nu -k} (\Gamma ; \Z )$ does not generate.
However, we have the following result of Ash and Rudolph:

\begin{theorem}\label{ashrudolph}
\cite{ash.rudolph} The restriction of $S_{0}\rightarrow H^{\nu
}(\Gamma; \Z )$ to the subgroup generated by the unimodular
symbols is surjective.
\end{theorem}

\begin{proof}
We present the proof of \cite{ash.rudolph}.  It suffices to show that
any modular symbol is equivalent mod $\partial S_{1}$ to a sum of
unimodular symbols.

Let $\vv = [v_{1},\dots ,v_{n}]$, and suppose that $\|\vv
\| > 1$.  Let $w\in \Z ^{n}$ be any point not in the lattice generated
by the $v_{i}$.  (Such a point exists since $\| \vv \| >1$.)  Let $\vv
_{i}$ be the modular symbol obtained by replacing $v_{i}$ with $w$ in
$\vv $.  Applying relation (4) from Theorem~\ref{msymb}, we have
\begin{equation}\label{fund.relation}
\vv = \sum (-1)^{i+1}\vv _{i}
\end{equation}   
in $S_{0}/\partial S_{1}$.  We claim $w$ can be modified so that
$0\leq \|\vv _{i}\| < \|\vv \|$, and at least one $\vv _{i}$ satisfies
$\|\vv _{i}\| \not = 0$.  This proves the theorem, because after
repeating the argument finitely many times, we can write $\vv $ as a
sum of unimodular symbols.

To prove the claim, write $w = \sum q_{i}v_{i}$, where $q_{i}\in \Q $.  We
have $\|\vv _{i}\| = |q_{i}|\|\vv \|$.  If we modify $w$ by
subtracting integral multiples of the $v_{i}$, we can ensure $0\leq
|q_{i}|<1$.  Furthermore, at least one $q_{i}\not =0$ since $w$ was
originally chosen not to lie in the lattice generated by the $v_{i}$.
\end{proof}

\subsection{}
Given a modular symbol $\vv $, the set of \emph{candidates} of $\vv $
is the set
\[
\Cand \vv = \Bigl\{w\in \Z ^{n}\Bigm|\text{$w\not =0$ and $w=\sum
q_{i}v_{i}$, where $0\leq
|q_{i}|<1$} \Bigr\}.
\]
The set $\Cand \vv $ contains exactly the points that may be used to
construct the homology \eqref{fund.relation} so that the resulting
modular symbols are closer to unimodularity.

For application of Theorem~\ref{ashrudolph} to Hecke
eigenvalue computations, we need to construct a candidate for
any $\vv $ with $\|\vv \|>1$.  We now discuss two conjectural
ways to do this.  These are useful for three reasons:
\begin{enumerate}
\item The conjectures will play an important role in
our algorithm to compute the Hecke
action on $H^{\nu -1} (\Gamma ;\Z )$.
\item The candidates produced by these methods are efficient in
practice,
in the sense that $\|\vv _{i}\|$ from \eqref{fund.relation} will be
much smaller than $\|\vv \|$. 
\item Conjecture~\ref{lllconj} provides an explicit polynomial-time
implementation of the modular symbol algorithm.
\end{enumerate}

Write $\vv =[v_{1},\dots ,v_{n}]$, and let $b (\vv )$ be the point
$\sum v_{i}v_{i}^{t}$.  One can show $b (\vv )\in C$ since
$\|\vv \| \not = 0$.  Recall that if $\sigma \in \V $, then $\Vertex
\sigma \subset \PPP $ is the set of primitive points corresponding to
the face of $\Pi $ that generates $\sigma $ (\S\ref{voronoi}).

\begin{conjecture}\label{vorconj}
Let $\vv $ be a modular symbol with $\|\vv \|>1$.  Let $\sigma \in \V$
be a top-dimensional cone containing $b(\vv)$.  Then
\[
\Cand \vv \cap \Vertex \sigma \not = \varnothing. 
\]
\end{conjecture}

\begin{remark}\label{vor.reduction.alg}
The cone $\sigma $ can be computed using the \emph{\Vor \ reduction
algorithm}~\cite[\S27ff]{voronoi1}.
\end{remark}

\subsection{}
Although geometrically attractive, the use of Conjecture~\ref{vorconj}
in practice suffers from two disadvantages.  First, to the best of our
knowledge, the complexity of the \Vor \ reduction algorithm is
unknown.  Second, the structure of $\Pi $ is difficult to
determine.\footnote{However, for $n\leq 4$, the structure of $\Pi $ is
well understood.  An elegant technique to index the faces using
configurations in projective space (in the sense of \cite{geom.imag})
can be found in \cite{mcc.art}.  To the best of our knowledge, the
complete structure of $\Pi $ is unknown for any other $n$, although
much is known for $5\leq n\leq 8$ (cf. \cite{conway.sloane.iii} and
the references there).}  An alternative uses $LLL$-reduction, which we
now recall.

\begin{definition}\label{lll.red}
\cite[Ch. 2.6]{henri.cohen}
Let $B=\{b_{1},\ldots,b_{n} \}$ be an ordered basis of $\R ^{n}$, and
let 
$B^{*}=\{b_{1}^{*},\ldots,b_{n}^{*} \}$ be the orthogonal (not
orthonormal) basis obtained from $B$ using the
Gram-Schmidt process.  Let  
$$\mu _{i,j} = ( b_{i}\cdot b_{j}^{*})/(b_{j}^{*}\cdot b_{j}^{*}),\quad 
\hbox{where $1\leq j<i\leq n$}.$$
Then $B$ is {\em $LLL$-reduced} if the
following inequalities hold:
\begin{enumerate}
\item $|\mu _{i,j}|\leq 1/2$, for $1\leq j<i\leq n$.
\item $|b_{i}^{*}+\mu _{i,i-1}b_{i-1}^{*}|^{2}\geq (3/4) |b_{i-1}^{*}|^{2}$.
\end{enumerate}
Furthermore, a quadratic form is said to be $LLL$-reduced if it is the Gram
matrix of an $LLL$-reduced basis.
\end{definition}

We emphasize that the basis $B$ in Definition~\ref{lll.red} is ordered.
Changing the order of $B$ changes $B^{*}$, which affects the
conditions of the definition.

\begin{conjecture}\label{lllconj}
Let $\vv $ be a modular symbol with $\|\vv \|>1$, and suppose that
$b(\vv )$ is an $LLL$-reduced quadratic form.  Let $E$ be the standard
basis for $\Z ^{n}$.  Then
\[
\Cand \vv \cap E \not =\varnothing .
\]
\end{conjecture}

\begin{remark}
To apply Conjecture~\ref{lllconj} in practice, one finds a matrix
$\gamma \in GL_{n} (\Z )$ such that $b (\gamma \cdot \vv) $ is
$LLL$-reduced, and then a candidate for $\vv $ will be in $\gamma
^{-1} E$.
\end{remark}

\subsection{}
We can prove the conjectures in some cases.  We begin by describing a
geometric interpretation of what it means for $w\in E$ to be a
candidate for $\vv $.

Let $\vv =[v_{1},\dots ,v_{n}]$ be a modular symbol, and fix an
ordering of the $v_{i}$.  Let $A$ be the matrix with columns $v_{i}$,
and let $B = \{b_{1},\ldots,b_{n} \}$ be the
basis made up of the rows of $A$.  Then one easily checks
that the quadratic form $b (\vv )$ is the Gram matrix of $B$.

\begin{lemma}\label{geo.interp}
Let $w = e_{k}\in E$.  For $1\leq i\leq n$, let $\vv _{i}$ be the
modular symbol constructed from $\vv $ and $w$ as in
\eqref{fund.relation}.  Also for $1\leq i\leq n$, let $B_{i}\subset \R
^{n-1}$ be the set of $(n-1)$ vectors obtained by projecting
$B\smallsetminus\{b_{k} \}$ into $P_{i}$, where $P_{i}$ is the span of
$E\smallsetminus\{e_{i} \}$.  Then the following statements are
equivalent:
\begin{enumerate}
\item $\|\vv _{i}\| < \|\vv \|$ for $1\leq i\leq n$.
\item $\vol B_{i} < \vol B$ for $1\leq i\leq n$.
\end{enumerate}
Here the volume in $P_{i}$ is normalized so that the fundamental
domains of $\Z ^{n}\subset \R ^{n}$ and $\Z ^{n}\cap P_{i}$ each have
volume $1$.
\end{lemma}

\begin{proof}
We have $\|\vv \| = \vol B = |\det A|$.  Furthermore, after choosing
$e_{k}$, we observe that $\|\vv _{i}\|$ and $\vol B_{i}$ are the
absolute value of the determinant of the same $(n-1)\times (n-1)$
minor of $A$.
\end{proof}

\begin{lemma}\label{new.lemma}
Let $\vv $ and $B$ be as above, and assume $\|\vv \|>1$.  If
$|b_{n}^{*}| > 1$, then $e_{n}\in \Cand \vv $.  
\end{lemma}

\begin{proof}
First note that $\vol B = \prod |b_{i}^{*}|$, since $B^{*}$ is
orthogonal.  Since $|b_{n}^{*}| > 1$ and $\vol B = \|\vv \|>1$, we
have 
\[
\prod _{i<n} |b_{i}^{*}| < \vol B.
\]

Now let $B_{i}$ be the projection of $B\smallsetminus \{b_{n} \}$ into
the coordinate hyperplane $P_{i}$, as in Lemma \ref{geo.interp}.
Clearly $\vol B_{i} \leq \prod _{i<n} |b_{i}^{*}|$.  Hence by Lemma
\ref{geo.interp}, $\|\vv_{i}\| < \| \vv \|$, and $e_{n}\in \Cand \vv $.

\end{proof}

\begin{proposition}\label{orthogonal}
Suppose $\|\vv \|>1 $ and $b(\vv )$ is a diagonal quadratic form.
Then Conjectures~\ref{vorconj} and~\ref{lllconj} are true.
\end{proposition}

\begin{proof}
First we show that Conjecture~\ref{lllconj} is true.  Since $b (\vv )$
is a diagonal quadratic form, we have $B=B^{*}$, and the $\mu _{ij}$
from Definition~\ref{lll.red} vanish.  Thus $\vol B = \prod
|b_{i}|>1$, and $|b_{i}|\geq 1$ for all $i$ since $B$ is integral.

Assume first that $B$ satisfies $|b_{i}| \leq |b_{j}|$ for $i\leq j$.
This implies $|b_{n}|>1$, and by Lemma~\ref{new.lemma} we have
$e_{n}\in \Cand \vv $, and Conjecture~\ref{lllconj} is true.

Now drop the assumption that $B$ is ordered by increasing lengths.  We
can multiply $\vv $ by a permutation matrix $\gamma $ so that $B$
satisfies $|b_{i}| \leq |b_{j}|$ for $i\leq j$.  This means that
$\gamma ^{-1}e_{n}\in \Cand \vv $.  Since $\gamma ^{-1}e_{n}\in E$,
Conjecture~\ref{lllconj} follows.

Finally, in this case Conjecture~\ref{lllconj} implies
Conjecture~\ref{vorconj}.  Since $b (\vv )$ is diagonal, it lies in
the cone $\sigma $ spanned by $\{q (e)\mid e\in E \}$.  This cone is a
proper face of the cone $\Sigma $ from Theorem~\ref{vor.thm}, and
hence $b (\vv )\in \Sigma$.  Since $E\subset \Vertex \Sigma$, the
result follows.
\end{proof}

Using standard estimates on $B$ and $B^{*}$, we can find a
lower bound on $\| \vv \| $ so that Conjecture~\ref{lllconj} is true.

\begin{proposition}\label{crude.bound}
Suppose that $\| \vv \| > 2^{n (n-1)/2}$.  Then Conjecture~\ref{lllconj}
is true.
\end{proposition}

\begin{proof}
We show that $\| \vv \| > 2^{n (n-1)/2}$ guarantees
$|b_{n}^{*}|>1$, which by Lemma~\ref{new.lemma} implies $e_{n}\in \Cand \vv $.
According to \cite[Theorem 2.6.2]{henri.cohen}, $B$ satisfies 
\[
\prod _{j} |b_{j}| \geq \| \vv \|  
\]
and 
\[
|b_{j}|\leq 2^{(n-1)/2}|b^{*}_{n}|,\quad \text{for $j=1,\ldots,n$.} 
\]
Hence
\[
2^{n(n-1)/2}|b^{*}_{n}|\geq\prod _{j} |b_{j}|\geq \| \vv \|.
\]
Solving for $|b^{*}_{n}|$, we see $\| \vv \| >
2^{n (n-1)/2}$ ensures $|b^{*}_{n}|>1$,
which proves the claim.
\end{proof}

\begin{theorem}\label{first.theorem}
Conjecture~\ref{vorconj} is true for $n=2$ and $3$.
\end{theorem}

\begin{proof}
We use Lemma~\ref{new.lemma} and direct investigation of the
reduction domains.  First we recall some facts about reduction theory
in these dimensions.  For convenience we use $GL_{n} (\Z )$ instead of
$SL_{n} (\Z )$.

For $n\leq 3$ the cone $\Sigma $ from Theorem~\ref{vor.thm} is the
only top-dimensional \Vor \ cone modulo $GL_{n} (\Z )$.  According to
\cite{conway.sloane.vi}, $b (\vv )\in \sigma $ if and only if $B$ is
an {\em obtuse superbase}.  By definition, this means the following.
Let $b_{0} = -\sum b_{i}$, and let $\bar B = B \cup \{b_{0} \}$.  Then
$\bar B$ satisfies
\[
b_{i}\cdot b_{j} \leq 0 \quad \text{for $0\leq i<j\leq n$}. 
\]

The set $\Sigma \cap C$ is not a fundamental domain for $GL_{n} (\Z )$
acting on $C$.  In fact, the stabilizer $\Gamma (\Sigma )\subset GL_{n}
(\Z )$ is a finite group, which for $n=2$ (respectively $3$) has order
$6$ (resp. $24$).  By placing additional conditions on the
basis $B$, we can describe a fundamental domain $T$ for $\Gamma (\Sigma )$
acting on $\Sigma $.

First we consider the case $n=2$.  The cone $\Sigma $ is a
$3$-dimensional cone inside the cone $\bar C$, and is spanned by $q
(e_{1})$, $q (e_{2})$, and $q (e_{1}-e_{2})$.  Figure \ref{red2-1.fig}
shows a $2$-dimensional affine slice of $\bar C$, with $\Sigma $
divided into fundamental domains for $\Gamma (\Sigma )$.  The shaded
region $T$ is half of the classical fundamental domain for $SL_{2} (\Z
)$ acting on $C$.

\begin{figure}[ht]
\psfrag{in}{$\infty $}
\psfrag{-1}{$-1$}
\psfrag{0}{$0$}
\psfrag{T}{$T$}
\begin{center}
\includegraphics[scale=.5]{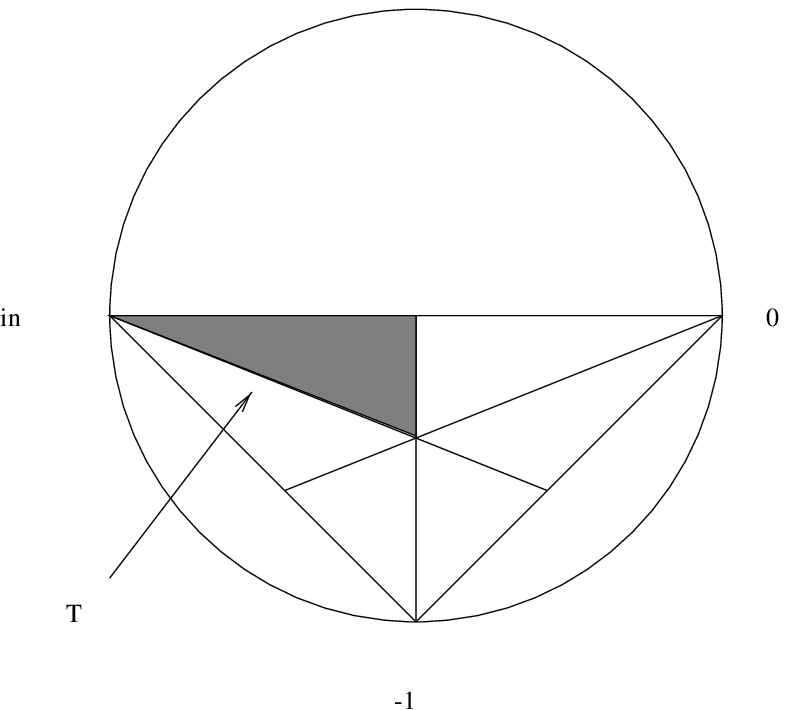}
\end{center}
\caption{\label{red2-1.fig}}
\end{figure}

Now we claim that if $b (\vv )\in T$ and $\|\vv \|\geq 2$, then
$e_{2}\in \Cand (\vv )$.  This implies the theorem for $n=2$,
because multiplying by elements of $\Gamma (\Sigma )$ stabilizes
$\Vertex \Sigma $.  

To prove the claim, we present another way to picture bases in the
region $T$.  If $B= (b_{1},b_{2})$, then $b (\vv )\in T$ if and only
if $B$ appears as in Figure \ref{red2-2.fig}.  In this figure, we have fixed
$b_{1}$, and $b_{2}$ must be in the infinite shaded region $S$ that lies above
the semicircle of radius $|b_{1}|$.  Points in $S$ correspond to ways
to complete $b_{1}$ to an obtuse superbase satisfying the additional
inequalities $|b_0|\geq |b_{2}|\geq |b_{1}|.$

\begin{figure}[ht]
\psfrag{S}{$S$}
\begin{center}
\includegraphics[scale=.5]{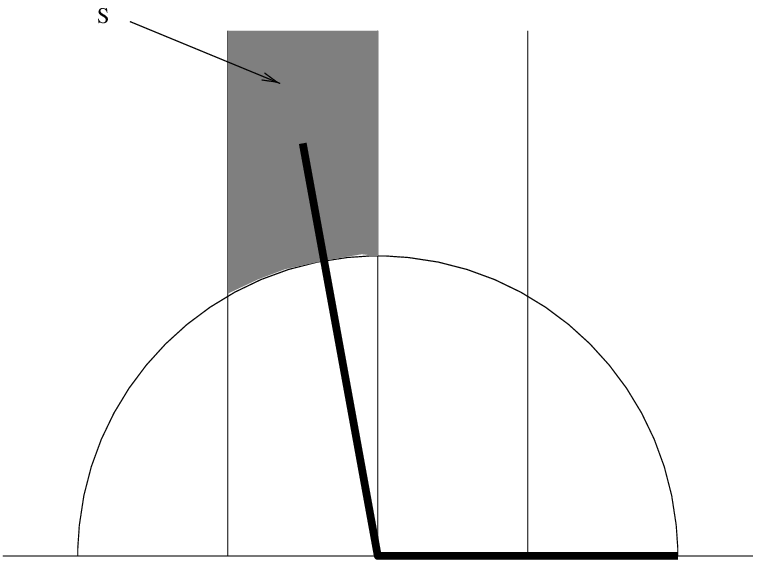}
\end{center}
\caption{\label{red2-2.fig}}
\end{figure}

Now consider the orthogonal basis $B^{*}$ constructed from $B$.  We
have $b_{1}=b_{1}^{*}$.  It
is easy to compute that $|b^{*}_{2}| \geq \sqrt{3}|b_{1}|/2$ for all
$b_{2}\in S$, and that the minimum occurs when $b_{2}$ is at the lower
left corner of $S$.  Hence if $|b_{1}|\geq \sqrt{2}$, we have
$|b^{*}_{2}|>1$, and by Lemma~\ref{new.lemma} we have $e_{2}\in \Cand
\vv $.

Since $B$ is integral, the remaining possibility is $|b_{1}| = 1$.
However, this implies that $b_{2}$ lies along the right edge of $S$,
and hence $b_{2}^{*}=b_{2}$.  If $|b_{2}|=1$, then $\|\vv \|=1$.  Thus
$|b_{2}|>1$ , and again $e_{2}\in \Cand \vv $.  This proves the
theorem for $n=2$.

The argument for $n=3$ is similar, although the reduction domain is
more complicated.  Now $\Sigma $ is $5$-dimensional, and the
fundamental domain $T$ can be described as follows.  As before, fix
$b_{1}$, and take $b_{2}$ to lie in the $2$-dimensional region $S$
from the $n=2$ case.  Together $b_{1}$ and $b_{2}$ determine the
Dirichlet-\Vor \ domain $pqrstu$ (see Figure~\ref{red3-1.fig}).  Let
$Z$ be the intersection of $pqrstu$ with 
\[
\{ x = \lambda _{1}b_{1} + \lambda _{2}b_{2}\mid x \cdot b_{1} \leq 0,
x\cdot b_{2}^{*}\leq 0\}.
\]
Then if $b (\vv )\in T$, the point $b_{3}$ must lie in the
$3$-dimensional region consisting of the points on or outside the
hemisphere of radius $|b_{2}|$ that project to $Z$.
Figure~\ref{red3-1.fig} shows the basis $B$, and Figure
\ref{red3-2.fig} shows $Z$ for different choices of $b_{2}$.
Altogether $T$ is a $5$-dimensional family of obtuse superbases that
can be described by additional inequalities similar to those for
$n=2$.

\begin{figure}[ht]
\psfrag{P}{$p$}
\psfrag{Q}{$q$}
\psfrag{R}{$r$}
\psfrag{S}{$s$}
\psfrag{T}{$t$}
\psfrag{U}{$u$}
\psfrag{V}{$V$}
\psfrag{b1}{$b_{1}$}
\psfrag{b2}{$b_{2}$}
\psfrag{b3}{$b_{3}$}
\begin{center}
\includegraphics[scale=.5]{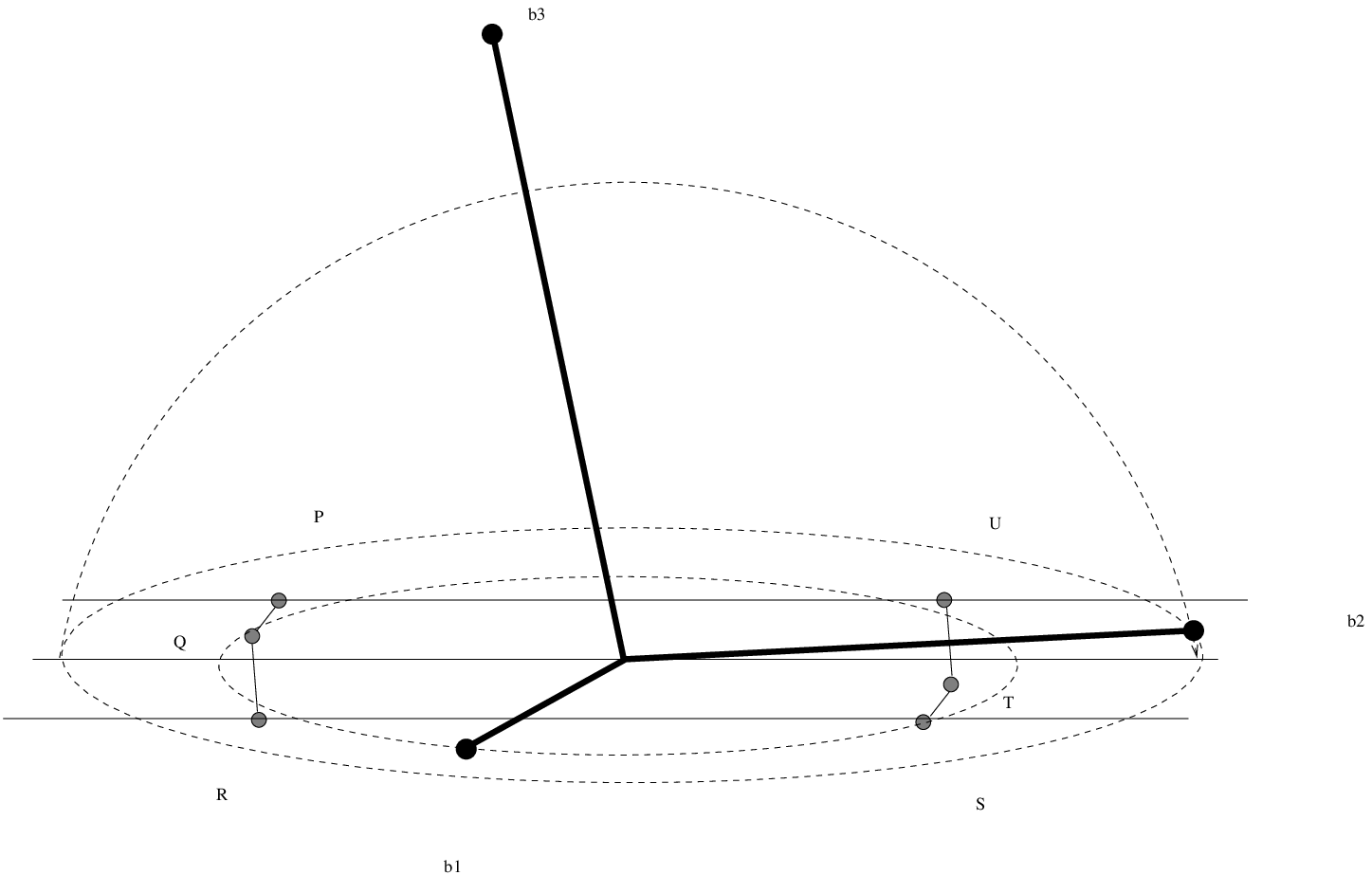}
\end{center}
\caption{\label{red3-1.fig}}
\end{figure}

\begin{figure}[ht]
\psfrag{Z}{$Z$}
\begin{center}
\includegraphics[scale=.5]{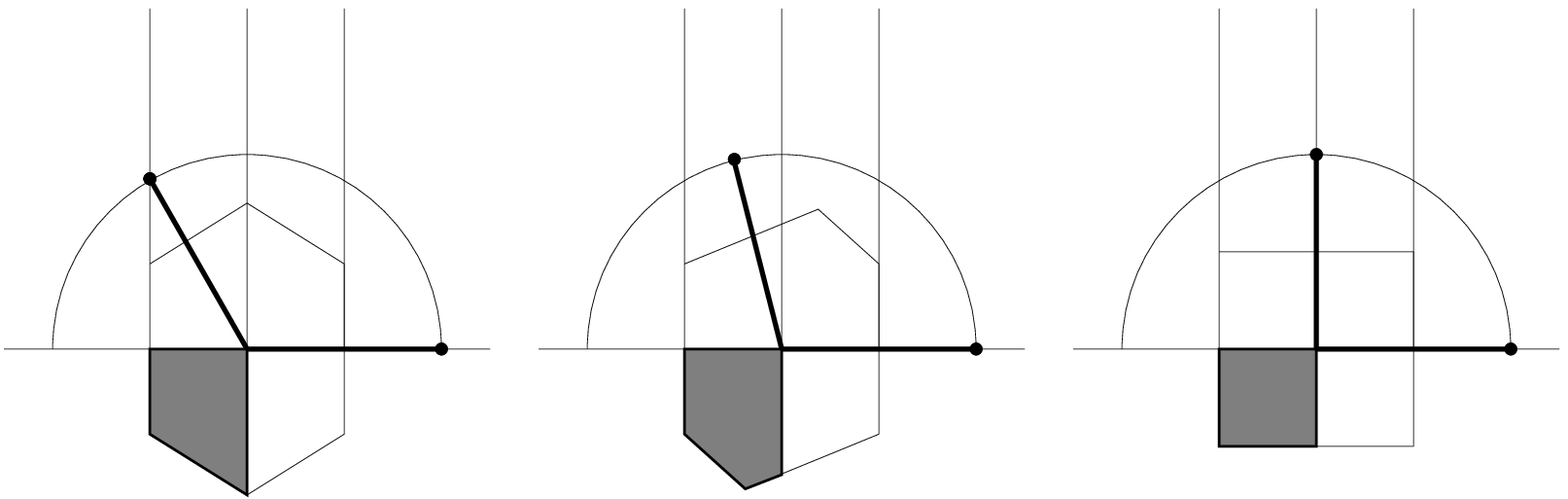}
\end{center}
\caption{\label{red3-2.fig}}
\end{figure}

We want to find conditions that imply $|b_{3}^{*}| > 1$, which will
imply $e_{3}\in \Cand \vv $.  Clearly the minimum value of
$|b_{3}^{*}|$ occurs when $|b_{1}|=|b_{2}| = |b_{3}|$.  Then for any
fixed $b_{2}$, the value of $|b_{3}^{*}|$ will be smallest when
$b_{3}$ projects to the vertices $a$ or $c$ of $Z$ shown in
Figure~\ref{red3-3.fig}.

So consider the set of bases satisfying 
\begin{enumerate}
\item $|b_{1}|=|b_{2}| = |b_{3}|$,
\item $0\geq b_{1}\cdot b_{2}\geq -|b_{1}|^{2}/2$, and 
\item $b_{3}$ projects to either $a$ or $c$ in
Figure~\ref{red3-3.fig}. 
\end{enumerate}

\begin{figure}[ht]
\begin{center}
\psfrag{a}{$a$}
\psfrag{c}{$c$}
\psfrag{a=c}{$a=c$}
\includegraphics[scale=.5]{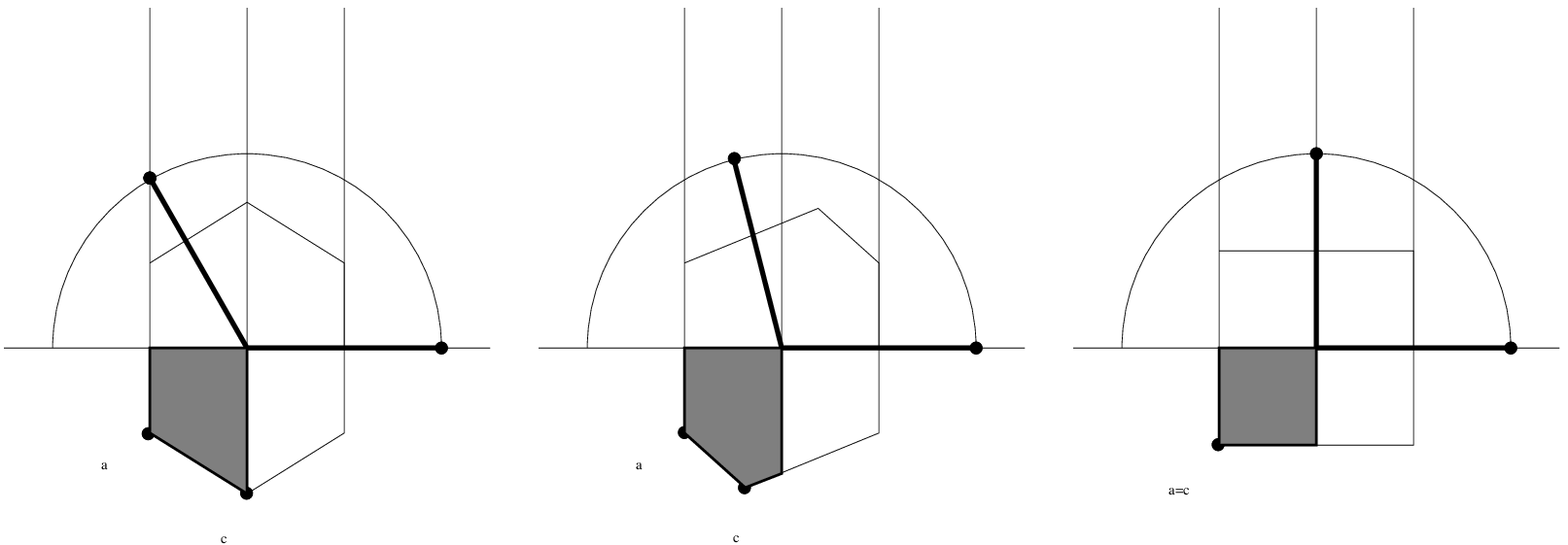}
\end{center}
\caption{\label{red3-3.fig}}
\end{figure}

It is not difficult to show that the minimal value of $|b^{*}_{3}|$ in
this family occurs when $a=c$, or when $b_{1}\cdot b_{2}=0$.  For this
basis, $|b_{3}^{*}|=|b_{1}|/\sqrt{2}$.  Hence if $|b_{1}|>\sqrt{2}$,
we have $e_{3}\in \Cand \vv $. 

The remaining cases are $|b_{1}| = 1$ or $\sqrt{2}$.  As for $n=2$, it
is straightforward, although tedious, to check that for any basis in
$T$ satisfying these conditions, we have either $\Cand \vv \cap
\Vertex \Sigma \not =\varnothing $ or $\|\vv \| = 1$.
\end{proof}

\begin{remark}\label{dee-n}
For $n=4$, there is only one other type of top-dimensional \Vor \ cone
mod $ GL_{4} (\Z )$, which corresponds to \Vor 's \emph{second perfect
form} \cite[\S34]{voronoi1}.  This cone corresponds to the lattice
$D_{4}$.  We are not aware of a useful characterization of the bases
appearing in this cone.
\end{remark}

%
%
\section{One-sharblies}\label{one.sharb}
In this section we describe our technique to compute the Hecke action
on $H^{\nu -1} (\Gamma ;\Z )$.

\subsection{}
Let $\xi = \sum n (\uu )\uu $ be a $k$-sharbly chain, where $n (\uu
)\in \Z $, and almost all $n (\uu )=0$.  Recall that a $k$-sharbly is
said to be reduced if and only if all its submodular symbols are
unimodular (Definition~\ref{reduced.sharbly}).

In general the reduced $k$-sharblies do not span $H^{\nu -k} (\Gamma
;\Z )$ (\S\ref{final.remarks}).  However, according to \cite{mcc.art},
$H^{\nu -1} (\Gamma ;\Z )$ is spanned by reduced $1$-sharblies if
$\Gamma \subset SL_{n} (\Z )$ and $n\leq 4$.  Hence to compute the
Hecke action on $H^{\nu -1} (\Gamma ;\Z )$ it suffices to describe an
algorithm that takes as input a $1$-sharbly cycle $\xi $ and produces
as output a cycle $\xi '$ satisfying:
\begin{enumerate}
\item [(a)] The classes of $\xi $ and $\xi '$ in $H^{\nu -1} (\Gamma ;\Z )$
are the same.
\item [(b)] $\|\xi' \| < \|\xi \|$ if $\|\xi \|>1$.
\end{enumerate}

We present an algorithm satisfying (a) in
Algorithm~\ref{hecke.algorithm}; in Conjecture~\ref{os.conj} we claim
the algorithm satisfies (b) for $n\leq 4$.  To simplify the
exposition, in \S\S\ref{sl2.start}--\ref{sl2.stop} we describe the
algorithm for $n=2$.  This case is arithmetically uninteresting---we
are describing how to compute the Hecke action on $H^{0} (\Gamma ; \Z
)$---but the geometry faithfully reflects the situation for all $n$.
We defer presentation for general $n$ to \S\ref{sln.start}.

\subsection{}\label{sl2.start}
Fix $n=2$, let $\xi \in S_{1}$ be a $1$-sharbly cycle mod $\Gamma$ for
some $\Gamma \subset SL_{2} (\Z )$, and suppose that $\xi $ is not
reduced.  We want to construct a cycle $\xi '$ homologous to $\xi $,
such that $\|\xi' \| < \|\xi \|$.  Since $\xi $ is not reduced, there
exist $\vv \in Z (\xi )$ with $\|\vv \|>1$.  Hence we want to perform
the modular symbol algorithm simultaneously over \emph{all} of
$\support \xi $ while constructing $\xi '$.  This leads to two
problems:
\begin{enumerate}
\item How should one choose candidates for the submodular symbols of
$\xi $?  Is the usual modular symbol algorithm sufficient for this?
\item Given $\xi $ and a collection of candidates for its submodular
symbols, how does one assemble the data into $\xi '$?
\end{enumerate}
Although these questions appear to be independent, they are in fact
coupled.  To answer the first, we claim that candidates should be
chosen using either Conjecture~\ref{vorconj} or~\ref{lllconj}; we
indicate why in \S\ref{interior.ms}.  We discuss the second in
\S\S\ref{gam.eqvariance}--\ref{endofdiscussion}.

\subsection{}\label{gam.eqvariance}
Suppose first that all $\vv \in Z (\xi )$ are nonunimodular.  We begin
by selecting candidates for each $\vv \in Z (\xi )$ using either
Conjecture~\ref{vorconj} or ~\ref{lllconj}, and we make these choices
$\Gamma $-equivariantly.  This means the following.  Suppose $\uu,\uu
'\in \support \xi $ and $\vv \in \support
(\partial \uu )$ and $\vv' \in \support (\partial \uu' )$ are modular
symbols such that $\vv = \gamma\cdot \vv ' $ for some $\gamma \in
\Gamma $.  Then we select $w\in \Cand \vv $ and $w'\in \Cand \vv '$
such that $w = \gamma\cdot w'$.

We can do this because if $\vv $ is a modular symbol and $w\in \Cand
\vv $, then $\gamma \cdot w \in \Cand (\gamma \cdot \vv )$ for any
$\gamma \in \Gamma $.  Since there are only finitely many $\Gamma
$-orbits in $Z (\xi )$, we can choose candidates $\Gamma
$-equivariantly by selecting them for some set of orbit
representatives.

It is important to note that $\Gamma $-equivariance is the only
``non-local'' criterion we use when selecting candidates.  In
particular, there is a priori no relationship among the $3$ candidates
chosen for any $\uu \in \support \xi $.

\subsection{}
Now we want to use the candidates and the $1$-sharblies in $\xi$ to
build $\xi '$.  Choose $\uu = [v_{1},v_{2},v_{3}]\in
\support \xi $, and denote the candidate for $[v_{i},v_{j}]$ by $w_{k}$,
where $\{i,j,k \} = \{1,2,3 \}$.  We use the $v_{i}$ and the
$w_{i}$ to build a $2$-sharbly chain $\eta (\uu )$ as follows.

\begin{figure}[ht]
\begin{center}
\psfrag{v1}{$v_{1}$}
\psfrag{v2}{$v_{2}$}
\psfrag{v3}{$v_{3}$}
\psfrag{w1}{$w_{1}$}
\psfrag{w2}{$w_{2}$}
\psfrag{w3}{$w_{3}$}
\includegraphics[scale=.5]{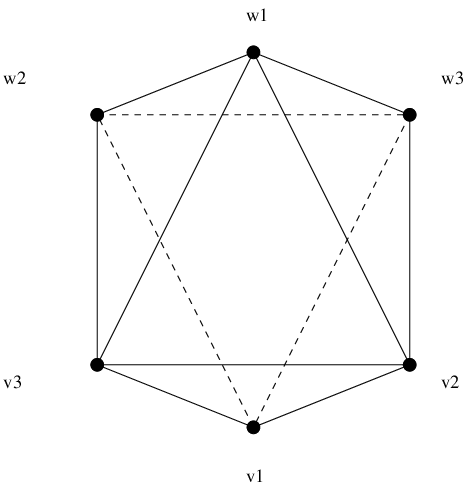}
\end{center}
\caption{\label{octa.fig}}
\end{figure}

\begin{figure}[ht]
\begin{center}
\psfrag{v1}{$v_{1}$}
\psfrag{v2}{$v_{2}$}
\psfrag{v3}{$v_{3}$}
\psfrag{w1}{$w_{1}$}
\psfrag{w2}{$w_{2}$}
\psfrag{w3}{$w_{3}$}
\includegraphics[scale=.5]{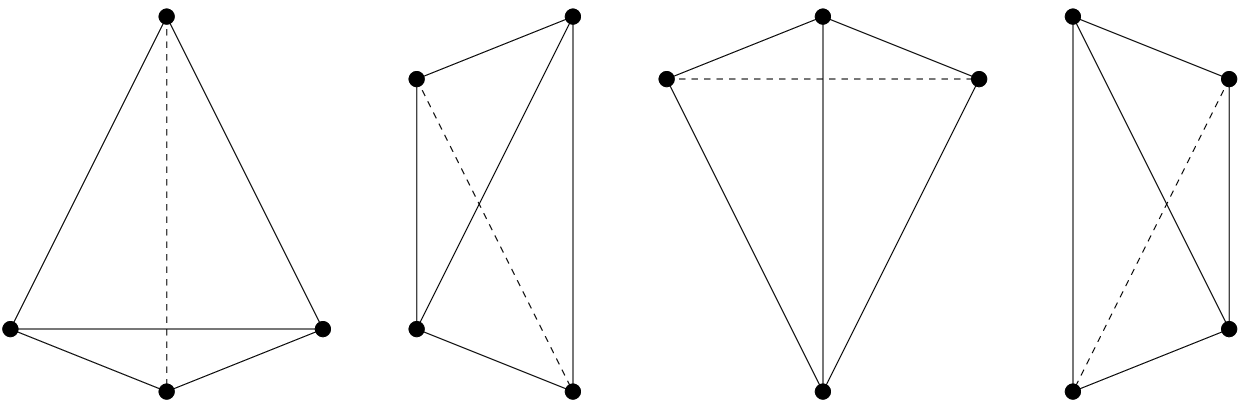}
\end{center}
\caption{\label{octa-hom.fig}}
\end{figure}

Let $P$ be an octahedron in $\R ^{3}$.  Label the vertices of $P$ with
the $v_{i}$ and $w_{i}$ such that the vertex labelled $v_{i}$ shares
no edge with the vertex labelled $w_{i}$ (Figure \ref{octa.fig}).  Now
subdivide $P$ into four tetrahedra without adding new vertices.  This
can be done by connecting two opposite vertices, say those with labels
$v_{1}$ and $w_{1}$, by a new edge (Figure \ref{octa-hom.fig}).

Now use the four tetrahedra to construct $\eta (\uu )$ as follows.
For each tetrahedron $T$, take the labels of four vertices and arrange
them into a quadruple.  If we orient $P$, then we can use the induced
orientation on $T$ to order the four primitive points.  In this way,
each $T$ determines a $2$-sharbly, and $\eta (\uu)$ is defined to be
the sum.  For example, if we use the decomposition in Figure
\ref{octa-hom.fig}, we have
\begin{equation}\label{eta}
\eta (\uu ) = [v_{1},v_{2},w_{3},w_{1}] + [v_{1},w_{3},w_{2},w_{1}] + [v_{1},w_{2},v_{3},w_{1}] +[v_{1},v_{3},v_{2},w_{1}].
\end{equation}

Now repeat this construction for all $\uu \in \support \xi $, and let 
$\eta =\sum n (\uu ) \eta (\uu )$.  Finally, let $\xi ' = \xi +\partial \eta $.

\subsection{}\label{endofdiscussion}
By construction, $\xi' $ is a cycle mod $\Gamma$ in the same
class as $\xi $.  We claim in addition that no submodular symbols from
$\xi $ appear in $\xi '$.  To see this, consider $\partial \eta (\uu
)$.  From \eqref{eta}, we have
\begin{multline}\label{bound}
\partial \eta (\uu ) = - [v_{1},v_{2},v_{3}] + [v_{1},v_{2},w_{3}] +
[v_{1},w_{2},v_{3}] + [w_{1},v_{2},v_{3}] \\
- [v_{1},w_{2},w_{3}] -
[w_{1},v_{2},w_{3}] - [w_{1},w_{2},v_{3}] + [w_{1},w_{2},w_{3}].
\end{multline}
Note that this is the boundary in $S_{*}$, not $\coinv{*}$.
Furthermore, it's easy to see that $\partial \eta (\uu )$ is
independent of which pair of opposite vertices of $P$ we connected to
define $\eta (\uu )$.  

From \eqref{bound}, we see that in $\xi +\partial \eta $, the
$1$-sharbly $-[v_{1},v_{2},v_{3}]$ is canceled by $\uu \in \support
\xi $.  Consider the $1$-sharblies in \eqref{bound} of the form
$[v_{i},v_{j},w_{k}]$.  We claim these $1$-sharblies vanish in 
$\partial_{\Gamma } \eta $.

To see this, suppose that $\uu $,$\uu '\in \support \xi $,
and suppose $\vv = [v_{1},v_{2}]\in \support \partial \uu $ equals
$\gamma \cdot \vv '$ for some $\vv ' = [v_{1}',v_{2}']\in \support
\partial \uu '$.  Since the candidates were chosen $\Gamma
$-equivariantly, we have $w=\gamma \cdot w'$.  This means that the
$1$-sharbly $[v_{1},v_{2}, w]\in \partial \eta (\uu )$ will be
canceled mod $\Gamma $ by $[v_{1}',v_{2}', w']\in \partial \eta (\uu'
)$.  Hence, in passing from $\xi $ to $\xi '$, the effect in
$\coinv{*}$ is to
replace $\uu $ with \emph{four} $1$-sharblies in $\support \xi '$:
\begin{equation}\label{tfm}
[v_{1},v_{2},v_{3}]\longmapsto - [v_{1},w_{2},w_{3}] -
[w_{1},v_{2},w_{3}] - [w_{1},w_{2},v_{3}] + [w_{1},w_{2},w_{3}].
\end{equation}
Note that in \eqref{tfm}, there are no $1$-sharblies of the form
$[v_{i},v_{j},w_{k}]$.

\begin{remark}\label{imp1}
For implementation purposes, it is not necessary to explicitly
construct $\eta $.  Rather, one may work directly with \eqref{tfm}.
\end{remark}

\subsection{}\label{interior.ms}
Why do we expect $\xi ' $ to satisfy $\|\xi '\| < \|\xi \|$?  First of
all, in the right hand side of \eqref{tfm} there are no submodular
symbols of the form $[v_{i},v_{j}]$.  In fact, any submodular symbol
involving a point $v_{i}$ also includes a candidate used to reduce the
$[v_{i},v_{j}]$.

However, consider the submodular symbols in \eqref{tfm} of the
form $[w_{i},w_{j}]$.  Since there is no relationship among the
$w_{i}$, one has no reason to believe that these modular
symbols are closer to unimodularity than those in $\uu $.  Indeed, one
might expect that these modular symbols satisfy $\|[w_{i},w_{j}]\|\geq
\|\uu \|$.  This is the content of problem~2 from
\S\ref{sl2.start}.

We claim that---in practice---if one uses Conjecture~\ref{vorconj}
or~\ref{lllconj} to select candidates, then these new modular symbols
will be very close to unimodularity.  In fact, usually they are 
trivial or satisfy $\|[w_{i},w_{j}]\|=1$.  To us, it seems that
Conjectures~\ref{vorconj} and~\ref{lllconj} select candidates
``uniformly'' over $\support \xi $, although we will not attempt to
make this notion precise.

\begin{remark}\label{best.cand}
To ensure $\|\xi '\| < \|\xi \|$, one must also choose the \emph{best}
candidate offered by the conjectures in a suitable sense
(\S\ref{why.is.xiprime.better}).
\end{remark}
\subsection{}\label{sl2.stop}
In the previous discussion we assumed that no submodular symbols of
any $\uu \in \support \xi $ were unimodular.  Now we discuss what to
do if some are.  As before, pick candidates for the nonunimodular
symbols.  There are three cases to consider.  

First, all submodular symbols of $\uu $ may be unimodular.  In this
case there are no candidates, and \eqref{tfm} becomes
\begin{equation}\label{tfm2}
[v_{1},v_{2},v_{3}]\longmapsto [v_{1},v_{2},v_{3}].
\end{equation}

\begin{figure}[ht]
\begin{center}
\psfrag{v1}{$v_{1}$}
\psfrag{v2}{$v_{2}$}
\psfrag{v3}{$v_{3}$}
\psfrag{w3}{$w_{3}$}
\includegraphics[scale=.5]{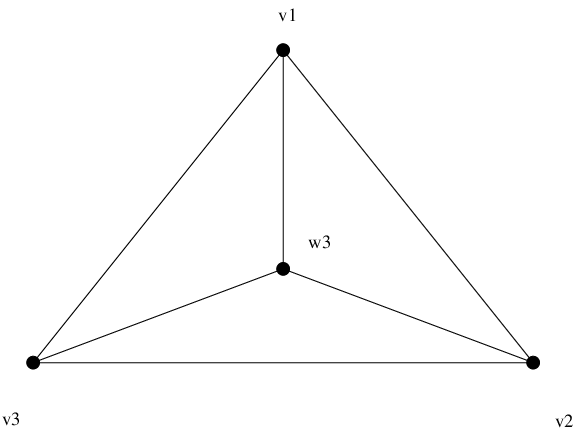}
\end{center}
\caption{\label{tetra.fig}}
\end{figure}

Second, one submodular symbol of $\uu $ may be nonunimodular, say the
symbol $[v_{1},v_{2}]$.  In this case we take $P$ to be a tetrahedron,
and $\eta (\uu ) = [v_{1},v_{2},v_{3},w_{3}]$
(Figure~\ref{tetra.fig}).  As before $[v_{1}, v_{2}, w_{3}]$ vanishes
in the boundary of $\eta $ mod $\Gamma $, and \eqref{tfm} becomes
\[
[v_{1},v_{2},v_{3}]\mapsto -[v_{1},v_{3},w_{3}] + [v_{2}, v_{3},
w_{3}].
\]

\begin{figure}[ht]
\begin{center}
\psfrag{v1}{$v_{1}$}
\psfrag{v2}{$v_{2}$}
\psfrag{v3}{$v_{3}$}
\psfrag{w2}{$w_{2}$}
\psfrag{w3}{$w_{3}$}
\includegraphics[scale=.5]{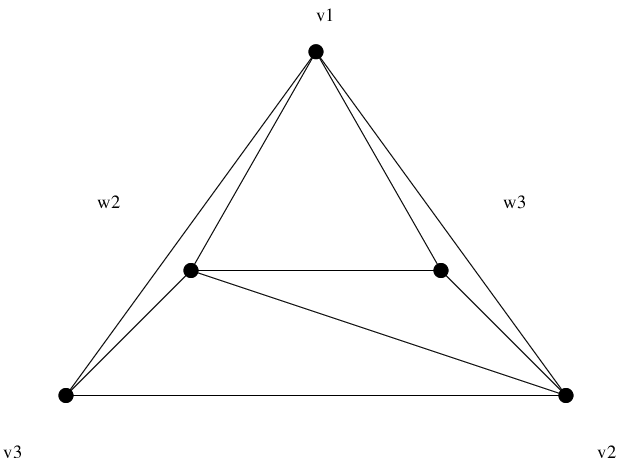}
\end{center}
\caption{\label{cone-on-square.fig}}
\end{figure}

Finally, two submodular symbols of $\uu $ may be nonunimodular, say
$[v_{1},v_{2}]$ and $[v_{1},v_{3}]$.  In this case we take $P$ to be
the cone on a square (Figure~\ref{cone-on-square.fig}).  To construct
$\eta (\uu )$ we must choose a decomposition of $P$ into tetrahedra.
Since $P$ has a non-simplicial face we must make a choice that affects
$\xi '$.  If we choose to subdivide $P$ by connecting the vertex
labelled $v_{2}$ with the vertex labelled $w_{2}$, we obtain
\[
[v_{1},v_{2},v_{3}]\longmapsto[v_{2},w_{2},w_{3}]+[v_{2},v_{3},w_{2}]+
[v_{1},v_{3},w_{2}].
\]

\subsection{}\label{sln.start}
We now describe the procedure for general $n$.  First we recall some
facts about convex polytopes.  Proofs can be found in \cite{ziegler}.

Let $P$ be a $d$-dimensional convex polytope embedded in $\R ^{d}$.
The \emph{facets} of a $d$-polytope $P$ are the faces of dimension
$(d-1)$.  The \emph{cone on $P$} is the polytope $cP$ constructed
as follows.  Choose a linear embedding $i\colon \R ^{d}\rightarrow \R
^{d+1}$ and let $x\in \R ^{d+1}\smallsetminus\R ^{d}$.  Then $cP$ is
the convex hull of $x$ and $i (P)$.  One can show that the
combinatorial type of $cP$ is independent of the choice of $x$ or $i$.
We also write $c^{0}P := P$ and $c^{k} (P) := c (c^{k-1}P)$.

Let $E$ be the standard basis of $\R ^{n}$.  Then the
\emph{$(n-1)$-simplex} $\Delta _{n-1}$ is the convex hull of $E$, and
the \emph{$n$-crosspolytope} $\beta _{n}$ is the convex hull of
$-E\cup E$.  Write $E=\{e_{i} \}$, and let $P (n,j)$ be the convex
hull of $E$ and the $j$ points $\{-e_{k}\mid 1\leq k\leq j\leq n \}$.

\begin{lemma}\label{it.cone}
The polytope $P (n,j)$ is isomorphic to the iterated cone
$c^{n-j}\beta _{j}$.
\end{lemma}

\begin{proof}
By definition, the convex hull of $A := \{\pm e_{k}\mid 1\leq k\leq
j\}$ is $\beta _{j}$.  The remaining vertices of $P (n,j)$ are the
points $B:= \{-e_{k}\mid j+1\leq k\leq n \}$.  Since $B$ is linearly
independent, and is also linearly independent of the linear span of
$A$, the lemma follows easily by induction.
\end{proof}

\begin{lemma}\label{subdivisions}
There are $j$ distinct subdivisions of $P (n,j)$ into simplices
without adding new vertices.
\end{lemma}

\begin{proof}
This follows immediately from Lemma~\ref{it.cone}.  Any such subdivision of
$\beta _{j}$ is formed by connecting one of the $j$ pairs of vertices
not already connected by an edge of $\beta _{j}$, and any such
subdivision of $c^{n-j}\beta _{j}$ is formed by subdividing $\beta
_{j}$ first. 
\end{proof}

\begin{algorithm}\label{hecke.algorithm}
Let $\Gamma $ be a torsion-free subgroup, and let $\xi  = \sum n (\uu
)\eta (\uu )$ be a $1$-sharbly
cycle mod $\Gamma $ representing a class in $H^{\nu -1} (\Gamma ;\Z
)$.  The output of this algorithm is a class $\xi '\in H^{\nu -1}
(\Gamma ;\Z )$.
\begin{enumerate}
\item [\textbf{A.}] \emph{Choose candidates.}  For each $\uu \in
\support \xi $, and for each $\vv \in \support \partial \uu $ with
$\|\vv \|>1$, choose a candidate $w (\vv )$.  Make these choices
$\Gamma $-equivariantly over all of $\support \xi $ as in
\S\ref{gam.eqvariance}.  For each $\uu \in \support \xi $, we let $C (\uu
)$ be the set $\{w (\vv )\mid \vv \in \support \partial \uu \}$.
\item [\textbf{B.}] \emph{Shift candidates.}  Choose $\uu=[v_{1},\dots
,v_{n+1}]\in \support \xi$, and set $j=\#C (\uu )$.  Apply relation
(1) from Definition \ref{sharbly.complex} so that the $j$ submodular
symbols
\[
\bigl\{\vv^{i} = [v_{1},\dots ,\hat v_{i},\dots ,v_{n+1}]\bigm|1\leq i\leq
j\bigr\}
\]
satisfy
$\|\vv^{i} \|>1$.  Write $w_{i}$ for $w (\vv^{i})$, and let $n' (\uu
)$ be the new coefficient of $\uu $ in $\xi $.
\item [\textbf{C.}] \emph{Construct $2$-sharblies.}  Let $P=P (n+1,j)$
be the polytope from Lemma~\ref{it.cone}, and choose a subdivision of
$P$ into simplices without adding new vertices as in Lemma
\ref{subdivisions}.  Orient $P$ so that the induced orientation on the
face spanned by $e_{1},\dots ,e_{n+1}$ is the opposite of the
orientation given by the standard ordering of $e_{1},\dots ,e_{n+1}$.
Via the correspondence
\begin{align*}
e_{k}&\longleftrightarrow v_{k}\quad \text{for $1\leq k\leq n+1$,}\\
-e_{k}&\longleftrightarrow w_{k}\quad \text{for $1\leq k\leq j$,}
\end{align*}
and the orientation of $P$, use the subdivision of $P$ to construct a
$2$-sharbly chain $\eta (\uu )$.
\item [\textbf{D.}] \emph{Continue.}  Complete steps \textbf{B} and
\textbf{C} for all $\uu \in \support \xi $.
\item [\textbf{E.}] \emph{Terminate.}  Set
\[
\eta = \sum _{\uu \in \support \xi} n' (\uu )\eta (\uu ),
\]
and define $\xi ' := \partial \eta +\xi $.
\end{enumerate}
\end{algorithm}

\subsection{}
Now we want to describe how $\xi '$ is related to $\xi$, and in
particular in what sense $\xi '$ is closer to unimodularity than $\xi
$.  Let $\uu \in \support \xi $, and let $\eta (\uu )$ be the
$2$-sharbly chain constructed above.  Define
\begin{align*}
\partial \eta _{\text{old}} (\uu ) &= - [v_{1},\dots ,v_{n+1}],\\
\partial \eta _{\text{side}} (\uu ) &= \sum _{k=1}^{j} [v_{1},\dots ,\hat{v_{k}},\dots ,v_{n+1},w_{k}],\\
\partial \eta _{\text{new}} (\uu ) &= \partial \eta (\uu )-\partial
\eta _{\text{old}} (\uu ) - \partial \eta _{\text{side}} (\uu ).
\end{align*}
Note that $\partial \eta _{\text{old}} (\uu ) $ and $\partial \eta
_{\text{side}} (\uu ) $ contain all the submodular symbols of $\uu $
that are nonunimodular.

\begin{theorem}\label{main.thm}
The cycle $\xi '$ constructed in Algorithm~\ref{hecke.algorithm} is
homologous to $\xi $.  If $\uu \in \support \xi$ and $\vv \in \support
\partial \uu $ with $\|\vv \|>1$, then $\vv $ does not appear as
submodular symbol of $\xi '$ in the following sense: 
\[
\xi ' = \sum _{\uu \in \support \xi } n' (\uu ) \partial \eta _{\text{new}} (\uu ).
\]
\end{theorem}

\begin{proof}
It is clear that $\xi '$ is homologous to $\xi $ mod $\Gamma $.  To
see the rest of the statement, first note that we have chosen
orientations so that 
\[
\xi  + \partial \eta = \sum _{\uu \in \support \xi}n' (\uu ) (\partial
\eta _{\text{side}} (\uu ) + \partial \eta _{\text{new}} (\uu )).
\]
Hence we must show 
\[
\partial \eta _{\text{side}} (\uu ) = 0 \mod \Gamma.
\]
We claim this follows since the candidates we chosen $\Gamma
$-equivariantly over all of $\support (\xi )$.  Indeed, any
$1$-sharbly in $\support (\partial \eta _{\text{side}} (\uu ) )$ is
built from a certain candidate and a $0$-sharbly in $\support
(\partial \xi )$.  An investigation of the orientations we chose in
the construction of $\partial \eta $ and the fact that $\partial
_{\Gamma } (\xi ) = 0$ show that $\partial \eta _{\text{side}} (\uu )
= 0 \mod \Gamma$.
\end{proof}

\subsection{}\label{why.is.xiprime.better}
To conclude this section we discuss conditions under which we expect
$\|\xi' \| < \|\xi \|$.  First we clarify Remark~\ref{best.cand}.

Let $\vv $ be a modular symbol, and let $w\in \Cand \vv $.  Let $\{\vv
_{i} \}$ be the modular symbols from \eqref{fund.relation} constructed
using $\vv $ and $w$.  Define an integer $\mu (w)$ by 
\[
\mu (w) = \Max_{i=1,\dots ,n} \left\{\|\vv_{i} \| \right\}.
\]

\begin{definition}\label{good.cand.S}
Let $S\subset \Cand \vv $.  Then $w\in S$ is a \emph{good candidate
from S} if 
\[
\mu (w) = \Min_{w'\in S} \left\{\mu (w') \right\}.
\]
Furthermore, we say that $w$ is a good candidate chosen using
Conjecture~\ref{vorconj} (respectively Conjecture~\ref{lllconj}) if
$w$ is a good candidate for the (conjecturally nonempty) intersection
indicated in Conjecture~\ref{vorconj}
(resp. Conjecture~\ref{lllconj}).
\end{definition}
Note that good candidates are not necessarily unique.

\begin{conjecture}\label{os.conj}
Suppose $n\leq 4$, and let $\xi $ and $\xi '$ be as in Algorithm
\ref{hecke.algorithm}.  Assume that $\|\xi \| > 1$.  Then if each $w
(\vv )$ from step \textbf{A} of Algorithm~\ref{hecke.algorithm} is a good
candidate chosen using Conjecture~\ref{vorconj} or~\ref{lllconj}, then 
$\|\xi' \| < \|\xi \|$.
\end{conjecture}

%
%
\section{Experiments}\label{experiments.section}
We conclude by describing experiments we performed to test
Conjectures~\ref{vorconj},~\ref{lllconj}, and~\ref{os.conj}.  These
experiments were performed at MIT and Columbia at various times from
1995 to 1998, on Sun (SunOS) and Intel (Linux) workstations.  We are
grateful to these departments for making this equipment and support
available.

Before we describe the experiments, we remark that all trials
completed successfully, and no counterexamples to the conjectures were
found.

\newcommand{\testentry}{\smallskip\noindent$\bullet\quad $}

\subsection{}\label{modular.symbol.tests}
The first experiments we performed addressed Conjectures~\ref{vorconj}
and~\ref{lllconj}.  Because of implementation difficulties mentioned
immediately after Remark~\ref{vor.reduction.alg}, we were only able to test
Conjecture~\ref{vorconj} in dimensions $\leq 4$.  However, we were
able to test Conjecture~\ref{lllconj} in dimensions $\leq 40$, thanks
to $LLL$-reduction code available in GP-Pari and LiDIA.

\testentry We began by testing finding candidates for random modular
symbols for $SL_{n} (\Z )$.  A random square integral matrix $m$ was
constructed with entries chosen some fixed range.  If $\det m\not =0$,
then we attempted to find a candidate for the modular symbol formed
from the columns of $m$.  We tried to test matrices with small
determinant, since for these modular symbols the set of candidates is
small.  
\begin{enumerate}
\item For $n=4$ we verified Conjecture~\ref{vorconj} on approximately
$20000$ matrices.
\item For $2\leq n\leq 20$, we verified Conjecture~\ref{lllconj} on
approximately $20000$ matrices from each dimension, and for $21\leq
n\leq 40$ we tested Conjecture~\ref{lllconj} on approximately $1000$
matrices from each dimension.  In these tests we rejected those
matrices whose determinants were outside the range specified by
Proposition~\ref{crude.bound}.
\end{enumerate}

\testentry Instead of random modular symbols, we tested coset
representatives of the double cosets in \eqref{diag.cosets} for
different dimensions and values of $p$ and $k$.  We used the standard
coset representatives found in \cite{krieg}. 
\begin{enumerate}
\item For $T_{p} (1,3)$, $T_{p} (2,3)$, $T_{p} (1,4)$, and $T_{p}
(3,4)$, we tested all primes $p\leq 97$ using both conjectures (again
discarding those outside the range of Proposition~\ref{crude.bound}).
\item For $T_{p} (2,4)$, we tested all primes $p\leq 67$ using both
conjectures.
\item For dimensions $5\leq n\leq 10$, we verified
Conjecture~\ref{lllconj} on representatives of $T_{p} (1,n)$ for $p=2,3$.
\end{enumerate}

\testentry Finally, we performed complete reduction of random modular
symbols.  In the previous experiments, we only verified that a
candidate for a given modular symbol could be found using our
conjectures.  In this case, we stored the resulting modular symbols on
a stack and iterated the process until all modular symbols were
unimodular.  Due to the large number of modular symbols produced, we
limited our tests of Conjecture~\ref{lllconj} to dimensions $\leq 10$,
and tested only medium-sized determinants, typically with absolute
value less than $20$.  We verified Conjecture~\ref{vorconj} on
approximately $2000$ modular symbols and Conjecture~\ref{lllconj} on
approximately $1000$ modular symbols from each dimension.

\subsection{}\label{lift.section}
To test Conjecture~\ref{os.conj}, we wanted to mimic the experiments
in \S\ref{modular.symbol.tests}.  This cannot be done naively for the
following reason.  A single modular symbol is automatically a cycle
mod $\Gamma $, but for a $1$-sharbly chain $\xi $ to be a cycle mod
$\Gamma $, nontrivial conditions must be met.  Furthermore,
Algorithm~\ref{hecke.algorithm} uses these conditions in an essential
way to decrease $\|\xi \|$.  

This dilemma has two resolutions.  Either we must test
Conjecture~\ref{os.conj} on cycles for specific $\Gamma \subset SL_{n}
(\Z )$, or we must design an implementation of
Algorithm~\ref{hecke.algorithm} that is ``local,'' i.e. operates on a
single $1$-sharbly at a time.  The first solution is not feasible if
one wishes to test many $1$-sharbly cycles, because such cycles are
very difficult to construct.  Hence we must take the second approach.

\begin{definition}\label{lifts}
Let $\uu$ be a basis element of $S_{k}$.  Then a \emph{lift} for $\uu
$ is an $n\times (n+k)$ integral matrix $M$ with primitive columns
such that $[M_{1},\dots ,M_{n+k}] = \uu $, where $M_{i}$ is the $i$th
column of $M$.
\end{definition}

Let $\xi $ be a $k$-sharbly cycle mod $\Gamma $.  We claim that $\xi $
may be encoded as a finite collection of $4$-tuples
$(\uu , n (\uu ), \{\vv \}, \{M (\vv ) \})$, where:
\begin{enumerate}
\item $\uu \in \support \xi $.
\item $n (\uu )\in \Z $.
\item $\{\vv  \} = \support \partial \uu $. 
\item $\{M (\vv ) \}$ is a set of lifts for $\{\vv \}$.  These lifts
are chosen to satisfy the following $\Gamma $-equivariance condition.
Suppose that for $\uu ,\uu '\in \support \xi $ we have $\vv \in
\support (\partial \uu ) $ and $\vv '\in \support (\partial \uu ')$
satisfying $\vv = \gamma \cdot \vv '$ for some $\gamma \in \Gamma $.
Then we require $M (\vv )= \gamma M (\vv ')$.
\end{enumerate}

Clearly any cycle can be represented by such data, although the
representation is far from unique.

\subsection{}\label{inherit}
Let $\psi = (\uu , n (\uu ), \{\vv \}, \{M (\vv ) \})$ be a $4$-tuple
that is part of a cycle $\xi $.  We claim that we can choose
candidates for $\{ \vv \}$ that will the equivariance condition in
\S\ref{gam.eqvariance} without knowing the rest of $\xi $.

To see this, recall that a square matrix $M = (M_{ij})$ with $\det M
\not =0$ is in \emph{Hermite normal form} if $M_{ij} = 0$ for $i < j$,
and $0\leq M_{ij} < M_{ii}$ for $i>j$.  Furthermore, if $\det M > 0$,
then $M_{ii}>0$.  It is standard that for any $M$, the orbit $GL_{n}
(\Z )\cdot M$ contains only one element in Hermite normal form
\cite[2.4.2]{henri.cohen}. 

Now to choose a candidate $w$ for $\vv \in \support (\partial \uu)$,
we compute the Hermite normal form $M_{0} (\vv )$ of $M (\vv )$ first,
and input $M_{0} (\vv )$ into Conjecture~\ref{vorconj}
or~\ref{lllconj} to compute $w$.  If $M (\vv ) = \gamma M (\vv ')$,
then $M_{0} (\vv ) = M_{0} (\vv ')$.  Hence by using lifts we guarantee
that candidate selection is $\Gamma $-equivariant, even though the
choices are made locally.

\subsection{}
This means that we can think of a random $4$-tuple $\psi $ as being a
piece of some unknown cycle $\xi $ mod $\Gamma $, and can test
Algorithm~\ref{hecke.algorithm} by trying to write $\psi $ as a
collection of reduced $4$-tuples.  To complete the discussion, we must
say how lifts are chosen for the submodular symbols of $\partial \eta
(\uu )$ that survive to $\xi '$.  

\begin{definition}\label{innandout}
Let $\uu =[v_{1},\dots ,v_{n+1}]$ be a $1$-sharbly, and let 
\[
\vv ^{i}=[v_{1},\dots ,\hat v_{i},\dots ,v_{n+1}],\quad \hbox{for
$1\leq i\leq n+1$}
\]
be the submodular
symbols in $\support \partial \uu$.  Suppose that $\|\vv ^{i}\|>1$ for
$1\leq i\leq j\leq n+1$, and let $W=\{w_{i}\mid 1\leq i\leq j \}$ be
the set of candidates.  Let $U$ be the set $\{v_{1},\dots ,v_{n+1} \}\cup W$.
Let $\vv =[u_{1},\dots ,u_{n}]$ be a modular symbol with $u_{i}\in U$.
\begin{enumerate}
\item The modular symbol $\vv $ is called an \emph{outer submodular
symbol of $\uu $} if exactly one $u_{i}\in W$.
\item The modular symbol $\vv $ is called an \emph{inner submodular
symbol of $\uu $} if two or more $u_{i}\in W$.
\end{enumerate}
\end{definition}

\begin{figure}[ht]
\begin{center}
\psfrag{v1}{$v_{1}$}
\psfrag{v2}{$v_{2}$}
\psfrag{v3}{$v_{3}$}
\psfrag{w1}{$w_{1}$}
\psfrag{w2}{$w_{2}$}
\psfrag{w3}{$w_{3}$}
\includegraphics[scale=.7]{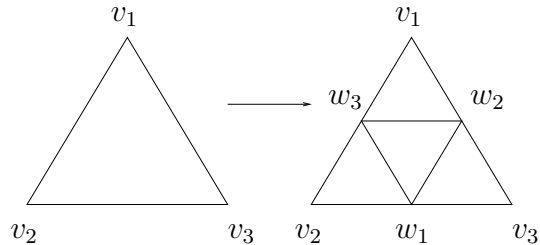}
\end{center}
\caption{$[w_{1},w_{2}]$ is inner, and $[v_{1},w_{3}]$ is
outer.\label{interior.fig}}
\end{figure}

Here is the meaning behind Definition~\ref{innandout}.  For
convenience suppose $n=2$ and $j=3$, and consider what happens when we
apply the algorithm to $\uu $.  We can think of $\uu $ as being a
triangle with vertices labelled by $v_{1}$, $v_{2}$, and $v_{3}$.
With this picture, to apply \eqref{tfm}, we can think of subdividing
the triangle into four new triangles, with the new vertices labelled
by the candidates $W$ (Figure \ref{interior.fig}).  

Now we discuss the relevance of inner and outer to our
implementation.  For an inner submodular symbol $\vv $, we can choose any
lift we like, as long as we choose the \emph{same}
lift for any other $1$-sharbly in \eqref{tfm} containing $\vv $.  If
$\vv $ is an outer submodular symbol, however, we must be more
careful.  In particular, consider Figure \ref{interior.fig}.  The lift
$M ([v_{1},v_{2}])$ was chosen using the $\Gamma $-action, and we must
choose $M ([v_{1},w_{3}])$ and $M ([v_{2},w_{3}])$ to reflect this.  

In practice, we can do the following.  If $\vv \in Z (\uu )$, then
each outer submodular symbol $\vv _{i}$ arising from $\vv $ is
obtained by replacing the $i$th primitive point of $\vv $ with $w$.
We construct $M (\vv _{i})$ by replacing the corresponding column of
$M (\vv )$ with $w$, and say that the lifts $\{M (\vv _{i}) \}$ are
\emph{inherited}. 

\begin{remark}\label{bad.reduction}
One might think that we could avoid computing Hermite normal forms and
could just apply Conjecture~\ref{vorconj} or~\ref{lllconj} directly.
However, this will not necessarily determine a unique representative
of the orbit $GL_{n} (\Z ) \cdot M (\vv )$, since this orbit may not
uniquely meet the \Vor \ and $LLL$ reduction domains.
\end{remark}

\subsection{}\label{one.sharbly.tests}
Now we describe the tests we performed to investigate
Conjecture~\ref{os.conj}.

\testentry We generated random $1$-sharblies $\xi $ with randomly
chosen lifts.  Using both modular symbol conjectures we constructed
candidates for $\xi $ and verified that $\|\xi '\|<\|\xi \|$.  Because
we only investigated dimensions $2$, $3$, and $4$, we were able to
test many $\xi $, approximately $10000$ per trial for $50$ trials.

\testentry We also tested all Hecke images within certain ranges
associated to certain ``standard'' reduced $1$-sharblies.  It is easy
to see that mod $SL_{n} (\Z )$, any reduced $1$-sharbly has the form
\begin{equation}\label{red.os.sln}
\left(\begin{array}{ccccc}
1&0&\dots &0&\varepsilon _{1}\\
0&1&\dots &0&\varepsilon _{2}\\
\vdots&&\ddots&\vdots&\vdots\\
0&0&\dots &1&\varepsilon _{n}
\end{array} \right),
\end{equation}
where the number of columns is $(n+1)$, and the last column is \[
(\varepsilon _{1},\dots ,\varepsilon _{n}) = (\underbrace{1,\dots
,1}_{k},0,\dots ,0), \quad \text{where $k=2,\dots ,n$.}
\]

Using these $1$-sharblies and randomly chosen lifts, we tested all
Hecke images within the following ranges:
\begin{enumerate}
\item For $T_{p} (1,3)$, $T_{p} (2,3)$, $T_{p} (1,4)$, and $T_{p}
(3,4)$, we tested all primes $p\leq 97$ using
Conjectures~\ref{vorconj} and ~\ref{lllconj}.
\item For $T_{p} (2,4)$, we tested all primes $p\leq 67$ using 
Conjectures~\ref{vorconj} and ~\ref{lllconj}.
\end{enumerate}
We repeated this experiment $10$ times to vary the lifts used.

\testentry We tested complete reduction of randomly chosen
$1$-sharblies with lifts.  At each step, the new $1$-sharblies
inherited lifts as described in \S\ref{inherit}.  This introduces the
possibility that for some initial choice of lifts, iteration of the
algorithm could fail to terminate.  However, this situation never
arose.  In 50 trials with approximately 10000 randomly chosen
$1$-sharblies, the complete reduction always terminated successfully.

\testentry After testing with random data, we computed the Hecke
action on cuspidal cycles occurring in $H^{2} (\Gamma _{0} (53); \Q )$,
where $\Gamma _{0} (53)\subset GL_{3} (\Z )$ is the subgroup of matrices
with bottom row equivalent to $(0,0,*) \mod 53$.  These cycles, or
rather their Lefschetz duals, were first discovered and investigated
in \cite{agg}.

We computed the characteristic polynomials of the Hecke operators
$T_{p} (1,3)$ for $p\leq 13$.  We found that these polynomials matched
those in \cite{agg}, which is consistent with the duality argument of
\cite[Theorem 3.1]{ash.tiep}.

\testentry Finally, in current work we are using the algorithm to
compute the Hecke action on $H^{5} (\Gamma_{0} (N); \Q )$, where
$\Gamma _{0} (N)\subset SL_{4} (\Z )$ is the subgroup of matrices with
bottom row congruent to $(0,0,0,*) \mod N$ \cite{agm}.  At the time of
this writing, we have completed computations for prime levels $N\leq
31$.  We have computed the characteristic polynomials for the Hecke
operators $T_{p} (k,4)$ for $1\leq k\leq 3$ and a range of $p$.  In
all cases the program wrote the Hecke image of a $1$-sharbly cycle as
sum of reduced $1$-sharbly cycles.

For these $GL_{3}$ and $SL_{4}$ tests, the author was helped and
encouraged enormously by Mark McConnell, who provided data for the
cycles generated by his program SHEAFHOM \cite{sheafhom}, and computed
the characteristic polynomials.

\subsection{}\label{final.remarks}
We conclude with a few remarks and open problems.

\testentry In general, if one wishes to implement the modular
symbol algorithm, Conjecture~\ref{lllconj} is much more efficient to
work with than Conjecture~\ref{vorconj}.  \Vor \ reduction is somewhat
difficult to program and requires a substantial amount of preliminary
computation.  On the other hand, high-quality computer code for
$LLL$-reduction is available from a variety of sources.

\testentry Algorithm~\ref{hecke.algorithm} can be adapted to work on
sharbly cycles $\xi \in S_{n+k}$ with $k>1$.  In particular, we can
describe the analogues of the polytopes $P (n,j)$ used in the
construction of $\xi '$: their facets involve
iterated cones on \emph{hypersimplices} \cite[Example 0.11]{ziegler}.  In
practice this is not useful for computing Hecke eigenvalues,
since we cannot expect in general that $\|\xi '\| < \|\xi \|$.

\testentry Throughout the description of Algorithm
\ref{hecke.algorithm}, we used the determinant as a measure of
``non-unimodularity'' of a $1$-sharbly.  Ultimately this approach
suffers from several shortcomings:

\begin{itemize}
\item [$\bullet$] For $\Gamma \subset SL_{n} (\Z )$ with $n\geq 4$, we
must use a nonreduced sharbly cycle to write a nontrivial element of
$H^{0} (\Gamma ;\Z )$.  
\item [$\bullet$] One wishes to compute Hecke eigenvalues in $H^{*}
(\Gamma ;\Z )$ for more exotic $\Gamma $.  For example, especially
interesting is $\Gamma \subset SL_{n} (\OK )$, where $\OK $ is the
ring of integers in a number field $K/\Q $.  If $\OK $ is not a
euclidean domain, then there is no obvious notion of a
\emph{primitive} vector.  One can still define the analogue of the
sharbly complex, and can use the determinant to define a $\Gamma
$-finite subset of sharblies \cite{bound}, but a practical modular
symbol algorithm is unknown in general.
\end{itemize}

A different approach is to use the relative position of a sharbly with
respect to $\Pi $ instead of the determinant.  This is carried out in
\cite{gunnells:modular} and \cite{sahc} for all arithmetic groups for
which $\Pi $ is available.  It would be nice to fuse the approach of
these articles and the approach described here.

\testentry If $\Gamma $ is not torsion-free, then our results hold if
we use cohomology with rational coefficients.  However, one can also
consider the equivariant cohomology $H_{\Gamma }^{*} (\Gamma ;\Z )$,
and can formulate conjectures about the arithmetic significance of
equivariant torsion classes \cite{ash.modp}.  Can
Algorithm~\ref{hecke.algorithm} be modified to compute the Hecke
action on $H^{\nu -1}_{\Gamma } (\Gamma ; \Z )$?

\testentry The modular symbol algorithm can be generalized to
$Sp_{2n}$ \cite{gunnells:symplectic}, and there is a cell complex that
can be used to compute $H_{*} (\Gamma )$, where $\Gamma \subset Sp_{4}
(\Z )$ \cite{macpherson.mcconnell:explicit}.  Is there a
``symplectic'' sharbly complex, and can an algorithm be devised to
compute Hecke eigenvalues on $H^{\nu -1} (\Gamma)$?

%
%
\bibliographystyle{amsplain}
\bibliography{experimental}

\end{document}